\documentclass[11pt, letterpaper, leqno]{amsart}
\usepackage{amsmath}
\usepackage{amsthm}
\usepackage{amsrefs}
\usepackage{lipsum}
\usepackage{amsopn}
\usepackage{amssymb}
\usepackage{fullpage}
\numberwithin{equation}{section}
\usepackage{float}
\usepackage{graphicx}
\usepackage{epstopdf}
\usepackage{titlesec}
\usepackage{fancyhdr}
\usepackage{fancyhdr}
\newtheorem{question}{Question}[section]

\newtheorem{lemma}{Lemma}[section]
\newtheorem{theorem}{Theorem}[section]

\newtheorem{definition}{Definition}[section]

\theoremstyle{definition}
\newtheorem{remark}{Remark}[section]

\DeclareMathOperator{\RE}{Re}
\DeclareMathOperator{\IM}{Im}
\DeclareMathOperator{\LOG}{Log}
\DeclareMathOperator{\Arg}{Arg}
\DeclareMathOperator{\adm}{adm}

\DeclareMathOperator{\arcsinh}{arcsinh}

\titleformat{\section}
  {\normalfont\scshape}{\thesection}{1em}{\centering}
\titleformat{\subsection}
{\normalfont\scshape}{\thesubsection}{1em}{\centering}

\begin{document}
\title{On a relation between harmonic measure and hyperbolic distance on planar domains}

\author{Christina Karafyllia}  
\address{Department of Mathematics, Aristotle University of Thessaloniki, 54124, Thessaloniki, Greece}
\email{karafyllc@math.auth.gr}   
\thanks{I thank Professor D. Betsakos, my thesis advisor, for his help and advice during the preparation of this work, and the Onassis Foundation for the scholarship I receive during my Ph.D. studies.}

\fancyhf{}
\renewcommand{\headrulewidth}{0pt}

\fancyhead[RO,LE]{\small \thepage}
\fancyhead[CE]{\small On a relation between harmonic measure and hyperbolic distance on planar domains }
\fancyhead[CO]{\small Christina Karafyllia} 
\fancyfoot[L,R,C]{}



\begin{abstract} Let $\psi $ be a conformal map of $\mathbb{D}$ onto an unbounded domain and, for $\alpha >0$, let ${F_\alpha }=\left\{ {z \in \mathbb{D}:\left| {\psi \left( z \right)} \right| = \alpha } \right\}$. If $\omega _\mathbb{D}\left( {0,{F_\alpha }} \right)$ denotes the harmonic measure at $0$ of $F_\alpha $ and $d_\mathbb{D} {\left( {0,{F_\alpha }} \right)}$ denotes the hyperbolic distance between $0$ and $F_\alpha$ in $\mathbb{D}$, then an application of the Beurling-Nevanlinna projection theorem implies that ${\omega _\mathbb{D}}\left( {0,{F_\alpha }} \right) \ge \frac{2}{\pi }{e^{ - {d_\mathbb{D}}\left( {0,{F_\alpha }} \right)}}$. Thus a natural question, first stated by P. Poggi-Corradini, is the following: Does there exist a positive constant $K$ such that for every $\alpha >0$, ${\omega _\mathbb{D}}\left( {0,{F_\alpha }} \right) \le K{e^{ - {d_\mathbb{D}}\left( {0,{F_\alpha }} \right)}}$? In general, we prove that the answer is negative by means of two different examples. However, under additional assumptions involving the number of components of $F_\alpha$ and the hyperbolic geometry of the domain $\psi \left( \mathbb{D} \right)$, we prove that the answer is positive.
\end{abstract}

\maketitle
\section{Introduction}\label{section1}

We will give an answer to a question of P. Poggi-Corradini (\cite[p. 36]{Co}) about an inequality relating harmonic measure and hyperbolic distance.  For a domain $D$, a point $z \in D$ and a Borel subset $E$ of $\overline D $, let ${\omega _D}\left( {z,E} \right)$ denote the harmonic measure at $z$ of $\overline E$ with respect to the component of $D \backslash {\overline E}$ containing $z$. The function ${\omega _D}\left( { \cdot ,E} \right)$ is exactly the solution of the generalized Dirichlet problem with boundary data $\varphi  = {1_E}$   (see \cite[ch. 3]{Ahl}, \cite[ch. 1]{Gar} and \cite[ch. 4]{Ra}). The hyperbolic distance between two points $z,w$ in the unit disk $\mathbb{D}$ (see \cite[ch. 1]{Ahl}, \cite[p. 11-28]{Bea}) is defined by 
\[{d_\mathbb{D}}\left( {z,w} \right) = \log \frac{{1 + \left| {\frac{{z - w}}{{1 - z\bar w}}} \right|}}{{1 - \left| {\frac{{z - w}}{{1 - z\bar w}}} \right|}}.\]
It is conformally invariant and thus it can be defined on any simply connected domain $D \ne \mathbb{C}$ as follows: If $f$ is a Riemann map of $\mathbb{D}$ onto $D$ and $z,w \in D$, then 
${d_D}\left( {z,w} \right) = {d_\mathbb{D}}\left( {{f^{ - 1}}\left( z \right),{f^{ - 1}}\left( w \right)} \right)$.
Also, for a set $E \subset D$, we define ${d_D}\left( {z,E} \right): = \inf \left\{ {{d_D}\left( {z,w} \right):w \in E} \right\}$.

The Hardy space with exponent $p$, $p>0$, and norm ${\left\|  \cdot  \right\|_p}$ (see \cite[p. 1-2]{Du}, \cite[p. 435-441]{Gar}) is defined to be
\[{H^p}\left( \mathbb{D} \right) = \left\{ {f \in H\left( \mathbb{D} \right):\left\| f \right\|_p^p = \mathop {\sup }\limits_{0 < r < 1} \int_0^{2\pi } {{{\left| {f\left( {r{e^{i\theta }}} \right)} \right|}^p}d\theta  <  + \infty } } \right\},\]
where $H\left( \mathbb{D} \right)$ denotes the family of all holomorphic functions on $\mathbb{D}$.  The fact that a function $f$ belongs to ${H^p}\left( \mathbb{D} \right)$ imposes a restriction on the growth of $f$ and this restriction is stronger as $p$ increases. If $\psi$ is a conformal map on $\mathbb{D}$, then $\psi \in {H^p}\left( \mathbb{D} \right)$ for all $p<1/2$ (\cite[p. 50]{Du}). Harmonic measure and hyperbolic distance are both conformally invariant and many Euclidean estimates are known for them. Thus, expressing the ${H^p}$-norms of a conformal map $\psi$ on $\mathbb{D}$ in terms of harmonic measure and hyperbolic distance, we are able to obtain information about the growth of the function by looking at the geometry of its image region $\psi \left( {\mathbb{D}} \right)$. Indeed, if $\psi $ is a conformal map on $\mathbb{D}$ and ${F_\alpha } = \left\{ {z \in \mathbb{D}:\left| {\psi \left( z \right)} \right| = \alpha } \right\}$ for $\alpha >0$, then (see \cite[p. 33]{Co})
\begin{equation}\label{isod}
\psi  \in {H^p}\left( \mathbb{D} \right) \Leftrightarrow \int_1^{ + \infty } {{\alpha ^{p - 1}}{\omega _{\mathbb{D}}}\left( {0,{F_\alpha }} \right)d\alpha }  <  + \infty.
\end{equation}

Now observe that if $E \subset {{\overline {\mathbb{D}}\backslash \left\{ 0 \right\}}}$, then ${\omega _\mathbb{D}}\left( {0,E} \right)$ and ${d _\mathbb{D}}\left( {0,E} \right)$ can be related by means of a special case of Beurling-Nevanlinna projection theorem (see \cite[p. 43-44]{Ahl}, \cite[p. 43]{Be}, \cite[p. 105]{Gar} and \cite[p. 120]{Ra}) which is stated as follows: Let $E \subset {{\overline {\mathbb{D}}\backslash \left\{ 0 \right\}}}$ be a closed, connected set intersecting the unit circle. Let ${E^ * } = \left\{ { - \left| z \right|:z \in E} \right\} = \left( { - 1,} \right.\left. { - {r_0}} \right]$, where ${r_0} = \min \left\{ {\left| z \right|:z \in E} \right\}$. Then,
	\[{\omega _{\mathbb{D}}}\left( {0,E} \right) \ge {\omega _{\mathbb{D} }}\left( {0,{E^ * }} \right) = \frac{2}{\pi }\arcsin \frac{{\left( {1 - {r_0}} \right)}}{{\left( {1 + {r_0}} \right)}}.\]

If $\psi $ is a conformal map of $\mathbb{D}$ onto an unbounded domain and
${F_\alpha } = \left\{ {z \in \mathbb{D}:\left| {\psi \left( z \right)} \right| = \alpha } \right\}$ for $\alpha >0$, then 
\[{d_\mathbb{D}}\left( {0,{F_\alpha }} \right) = \inf \left\{ {{d_\mathbb{D}}\left( {0,z} \right):z \in {F_\alpha }} \right\} = \log \frac{{1 + {r_0}}}{{1 - {r_0}}},\]
where ${r_0} = \min \left\{ {\left| z \right|:z \in {F_\alpha }} \right\}$. Thus, the Beurling-Nevanlinna projection theorem implies that

\begin{equation}\label{1.1}
{\omega _\mathbb{D}}\left( {0,{F_\alpha }} \right) \ge \frac{2}{\pi }\arcsin \frac{{\left( {1 - {r_0}} \right)}}{{\left( {1 + {r_0}} \right)}} = \frac{2}{\pi }\arcsin {e^{ - {d_\mathbb{D}}\left( {0,{F_\alpha }} \right)}} \ge \frac{2}{\pi }{e^{ - {d_\mathbb{D}}\left( {0,{F_\alpha }} \right)}}.
\end{equation}

Poggi-Corradini observed that, in general, the opposite inequality fails. But for a sector domain (\cite[p. 34-35]{Co}), 
\[{{\omega _{\mathbb{D}}}\left( {0,{F_\alpha }} \right) \le K{e^{ - {d_{\mathbb{D}}}\left( {0,{F_\alpha }} \right)}}}.\]
So, taking all these results into consideration, he set the following questions (\cite[p. 36]{Co}):  

\begin{question}\label{quest} Let $\psi $ be a conformal map  of $\mathbb{D}$ onto an unbounded domain and,  for $\alpha >0$, let ${F_\alpha } = \left\{ {z \in \mathbb{D}:\left| {\psi \left( z \right)} \right| = \alpha } \right\}$. Does there exist a positive constant $K$ such that for every $\alpha >0$,
	\[{\omega _\mathbb{D}}\left( {0,{F_\alpha }} \right) \le K{e^{ - {d_\mathbb{D}}\left( {0,{F_\alpha }} \right)}}?\]
\end{question}

\begin{question}\label{que} More generally, is it true that
	\[\psi  \in {H^p}\left( {\mathbb{D}} \right) \Leftrightarrow \int_1^{ + \infty } {{\alpha ^{p - 1}}{e^{ - {d_{\mathbb{D}}}\left( {0,{F_\alpha }} \right)}}d\alpha }  <  + \infty ?\]	
\end{question}

\begin{figure}[H]
	\begin{minipage}{0.4\textwidth}
		\begin{center}
			\includegraphics[scale=0.3]{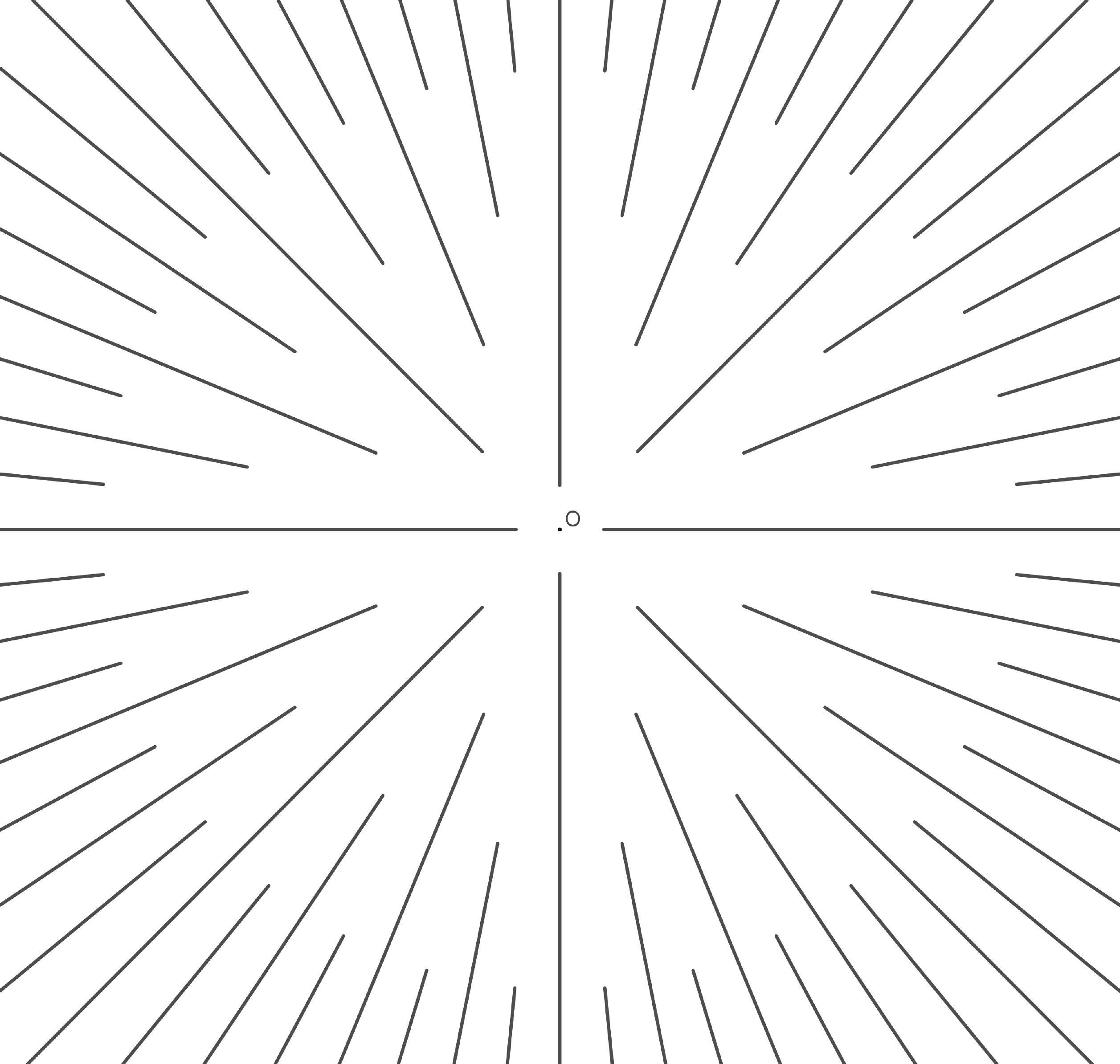}
			\vspace*{0.2cm}
			\caption{The simply connected domain $D$.}
			\label{t1}
		\end{center}
	\end{minipage}\hfill
	\begin{minipage}{0.6\textwidth}
		\begin{center}
			\includegraphics[scale=0.28]{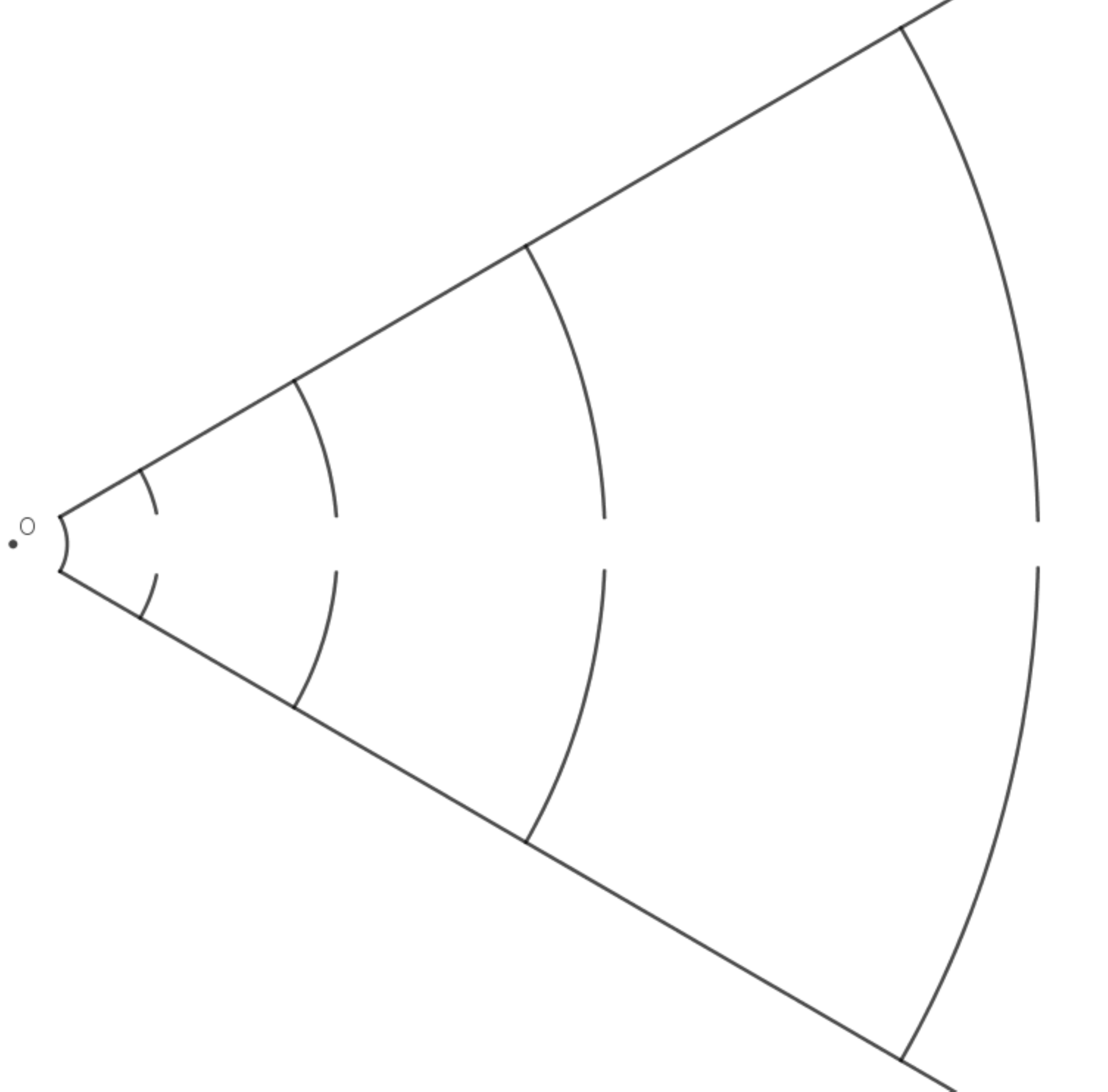}
			\caption{The simply connected domain $D'$.}
			\label{t2}
		\end{center}
	\end{minipage}
\end{figure}

In Section \ref{section4} we give a negative answer to the first question by mapping, through a conformal map $\psi$, $\mathbb{D}$ onto the simply connected domain $D$ of Fig. \ref{t1}. Its special feature is that as $\alpha  \to  + \infty $, the number of components of $\psi \left( {{F_\alpha }} \right)$ tends to infinity. This in conjunction with the fact that ${d_D}\left( {0,\psi \left( {{F_\alpha }} \right)} \right)$ is related to one component of $\psi \left( {{F_\alpha }} \right)$ whereas ${\omega _D}\left( {0,\psi \left( {{F_\alpha }} \right)} \right)$ is related to the whole $\psi \left( {{F_\alpha }} \right)$, made us believe that the choice of $D$ would give a negative answer to the Question \ref{quest} and so it did.

Consequently, a natural query would be whether the answer to the  Question \ref{quest} is positive in case the number of components of $\psi \left( {{F_\alpha }} \right)$ is bounded from above by a positive constant for every $\alpha>0$. However, in Section \ref{section5} we prove  by mapping $\mathbb{D}$ onto the simply connected domain $D'$ of Fig. \ref{t2}, that the answer is again negative. This is due to the fact that the hyperbolic distance between $\psi \left( {{F_\alpha }} \right)$ and the geodesic, $\psi \left( {{\Gamma _\alpha }} \right)$, joining the endpoints of $\psi \left( {{F_\alpha }} \right)$ in $D'$ tends to infinity, as $\alpha  \to  + \infty $. These results lead us to set sufficient conditions on the domain $\psi \left( {\mathbb{D}} \right)$ in order to give a positive answer to the Question \ref{quest}. 

In Section \ref{section2} we introduce some preliminary notions and results such as the domain decomposition method studied by N. Papamichael and N.S. Stylianopoulos (see \cite{Pa9}, \cite{Pap}, \cite{Pab}). In Section \ref{section3} we present some lemmas required for the proof of the theorem which is stated and proved in Section \ref{section4} and gives a negative answer to the Question \ref{quest} through the study of the domain of Fig. \ref{t1}. In Section \ref{section5} we consider a different kind of domain (see Fig. \ref{t2}) and prove that the answer is still negative. Having these results in mind we finally set the sufficient conditions to give a positive answer to the Question \ref{quest}. First note that if $\psi $ is a conformal map  of $\mathbb{D}$ onto an unbounded domain $D$ and ${F_\alpha } = \left\{ {z \in \mathbb{D}:\left| {\psi \left( z \right)} \right| = \alpha } \right\}$ for $\alpha>0$, then $\psi \left( {{F_\alpha }} \right)$ is a countable union of open arcs in $D$ which are the intersection of $D$ with the circle $\left\{ {z \in \mathbb{C}:\left| z \right| = \alpha } \right\}$ and have two distinct endpoints on $\partial D$. Thus, the preimage of every such arc is also an arc in $\mathbb{D}$ with two distinct endpoints on $\partial \mathbb{D}$ (see Proposition 2.14 \cite[p. 29]{Pom}). So, in Section \ref{section6} we prove the following theorem:

\begin{theorem}\label{kyrio} Let $\psi $ be a conformal map of $\mathbb{D}$ onto an unbounded simply connected domain $D$ and let ${F_\alpha } = \left\{ {z \in \mathbb{D}:\left| {\psi \left( z \right)} \right| = \alpha } \right\}$ for $\alpha >0$. If $N\left( \alpha  \right)$ denotes the number of components of $F_\alpha$ and ${F_\alpha }^j$ denotes each of these components for $j = 1,2, \ldots ,N\left( \alpha  \right)$, then we set $z_\alpha ^j,{z_\alpha ^j}'$ be the endpoints of ${F_\alpha }^j$ on $\partial \mathbb{D}$ and ${{\Gamma}_\alpha }^j$ be the geodesic joining $z_\alpha ^j$ to ${z_\alpha ^j}'$ in $\mathbb{D}$ for $j = 1,2, \ldots ,N\left( \alpha  \right)$. Suppose that the following conditions are satisfied:
\begin{enumerate}
\item There exists a positive constant $c_1$ such that $N\left( \alpha  \right) \le {c_1}$ for every $\alpha>0$.
\item There exists a positive constant $c_2$ such that $\psi \left( {{F_\alpha }^j} \right) \subset \left\{ {z \in D:{d_D}\left( {z,\psi \left( {{\Gamma _\alpha }^j} \right)} \right) < {c_2}} \right\}$ for every $\alpha>0$ and every $j \in \left\{ {1,2, \ldots ,N\left( \alpha  \right)} \right\}$.
\end{enumerate}
Then there exists a positive constant $K$ such that for every $\alpha>0$,
	\[{\omega _\mathbb{D}}\left( {0,{F_\alpha }} \right) \le K{e^{ - {d_\mathbb{D}}\left( {0,{F_\alpha }} \right)}}.\] 

\end{theorem}

\begin{remark} The direction ``$\Rightarrow$'' in Question \ref{que} is a simple consequence of (\ref{isod}) and (\ref{1.1}). So, the question actually concerns the other direction. The domain $D$ does not give an answer because by Theorem 4.1 in \cite[p. 239]{Ha}, which gives the Hardy number of an unbounded starlike region with respect to $z=0$, we derive that the corresponding Riemann map belongs in ${H^p}\left( \mathbb{D} \right)$ for every $p>0$. So, (\ref{isod}) implies that
\[\int_1^{ + \infty } {{\alpha ^{p - 1}}{\omega _{\mathbb{D}}}\left( {0,{F_\alpha }} \right)d\alpha }  <  + \infty\]
for every $p>0$.
\end{remark}

\begin{remark} The domain $D'$ does not give an answer to the Question \ref{que} because we have found that $\omega _{\mathbb{D}}\left( {0,{F_\alpha }} \right)$ decreases very rapidly so that
\[\int_1^{ + \infty } {{\alpha ^{p - 1}}{\omega _{\mathbb{D}}}\left( {0,{F_\alpha }} \right)d\alpha }  <  + \infty\]
for every $p>0$. This follows from calculations which we don' t present here.
\end{remark}

\section{Preliminary results and notations}\label{section2}

\subsection{Minda's reflection principle}

Concerning the hyperbolic metric we use the following theorem known as Minda's Reflection Principle \cite[p. 241]{Mi}. First, note that, if $\Gamma$ is a straight line (or circle), then $R$ is one of the half-planes (or the disk) determined by $\Gamma$ and ${\Omega ^*}$ is the reflection of a hyperbolic region $\Omega$ in $\Gamma$ .

\begin{theorem} \label{rp} Let $\Omega$ be a hyperbolic region in $\mathbb{C}$ and $\Gamma$ be a straight line or circle with $\Omega  \cap \Gamma  \ne \emptyset $. If $\Omega \backslash R \subset \Omega^*$, then 
	\[{\lambda _{{\Omega ^*}}}\left( z \right) \le {\lambda _\Omega }\left( z \right)\]
	for all $z \in \Omega \backslash \overline R $. Equality holds if and only if $\Omega$ is symmetric about $\Gamma$.
\end{theorem}

\subsection{Quasi-hyperbolic distance}

The hyperbolic distance between $z_1,z_2 \in D$ can be estimated by the quasi-hyperbolic distance, ${\delta _D}\left( {z_1,z_2} \right)$, which is defined by
\[{\delta _D}\left( {{z_1},{z_2}} \right) = \mathop {\inf }\limits_{\gamma :{z_1} \to {z_2}} \int_\gamma  {\frac{{\left| {dz} \right|}}{{d\left( {z,\partial D} \right)}}}, \]
where the infimum ranges over all the paths  connecting $z_1$ to $z_2$ in $D$ and $d\left( {z,\partial D} \right)$ denotes the Euclidean distance of $z$ from $\partial D$. Then it is proved that $\left( {{1 \mathord{\left/
			{\vphantom {1 2}} \right.
			\kern-\nulldelimiterspace} 2}} \right){\delta _D} \le {d_D} \le 2{\delta _D}$ (see \cite[p. 33-36]{Bea}, \cite[p. 8]{Co}).

\subsection{Extremal length}

Another  conformally invariant quantity which plays a central role in the proof of Section \ref{section4} is the extremal length. We present the definition and the properties we need as they are stated in \cite[ch. 4]{Ahl}, \cite[p. 361-385]{Be}, \cite[ch. 7]{Fu}, \cite[ch. 4]{Gar}, \cite[p. 88-100]{koo}  and \cite[ch. 2]{Oh}.

\begin{definition} Let $\left\{ C \right\}$ be a family of curves and $\rho \left( z \right) \ge 0$ be a measurable function defined in $\mathbb{C}$. We say $\rho \left( z \right)$ is admissible for $\left\{ C \right\}$ and denote by $\rho  \in \adm \left\{ C \right\}$, if for every rectifiable $C \in \left\{ C \right\}$, the integral $\int_C {\rho \left( z \right)\left| {dz} \right|} $ exists and $1 \le \int_C {\rho \left( z \right)\left| {dz} \right|}  \le  + \infty $. The extremal length of $\left\{ C \right\}$, $\lambda  \left\{ C \right\}$, is defined by
	\[\frac{1}{{\lambda \left\{ C \right\}}} = \mathop {\inf }\limits_{\rho  \in adm\left\{ C \right\}} \int \int {{\rho ^2}\left( z \right)dxdy}. \]
\end{definition}

Note that if all curves of $\left\{ C \right\}$ lie in a domain $D$, we may take ${\rho \left( z \right) = 0}$ outside $D$. The conformal invariance is an immediate consequence of the definition (see \cite[p. 90]{Fu}). As a typical example (see \cite[p. 366]{Be}, \cite[p. 131]{Gar}), we mention the case in which $R$ is a rectangle with sides of length $a$ and $b$ and $\left\{ C \right\}$ is the family of curves in $R$ joining the opposite sides of length $a$. Then $\lambda \left\{ C \right\} = \frac{b}{a}$. Next we state two basic properties of extremal length that we will need (see \cite[p. 54-55]{Ahl}, \cite[p. 363]{Be}, \cite[p. 91]{Fu}, \cite[p. 134-135]{Gar}, \cite[p. 79]{Oh}).

\begin{theorem}\label{el} If $\left\{ {C'} \right\} \subset \left\{ C \right\}$ or every $C' \in \left\{ {C'} \right\}$ contains a $C \in \left\{ {C} \right\}$, then $\lambda \left\{ C \right\} \le \lambda \left\{ {C'} \right\}.$
\end{theorem}

\begin{theorem}[The serial rule]\label{sr} Let $\left\{ {B_n} \right\}$ be mutually disjoint Borel sets and each $C_n \in \left\{ {C_n} \right\}$ be in $B_n$. If $ \left\{ {C} \right\}$ is a family of curves such that each $C$ contains at least one $C_n$ for every $n$, then
	\[\lambda \left\{ C \right\} \ge \sum\limits_n {\lambda \left\{ {{C_n}} \right\}}. \] 
\end{theorem} 

Sometimes it is more convenient to use the more special notion of extremal distance. Let $D$ be a plane domain and $E_1,E_2$ be two disjoint closed sets on $\partial D$. If $\left\{ {C} \right\}$ is the family of curves in $D$ joining $E_1$ to $E_2$, then the extremal length $\lambda_D \left\{ C \right\}$ is called the extremal distance between $E_1$ and $E_2$ with respect to $D$ and is denoted by ${\lambda _D}\left( {{E_1},{E_2}} \right)$.

\subsection{Domain decomposition method}

In case of quadrilaterals, the opposite inequality in the serial rule has been studied by N. Papamichael and N.S. Stylianopoulos by means of a domain decomposition method for approximating the conformal modules of long quadrilaterals (see \cite{Pa9}, \cite{Pap}, \cite{Pab}). Before stating the theorem we need, we present the required notation. 

\begin{figure}[H]
	\begin{center}
		\includegraphics[scale=0.5]{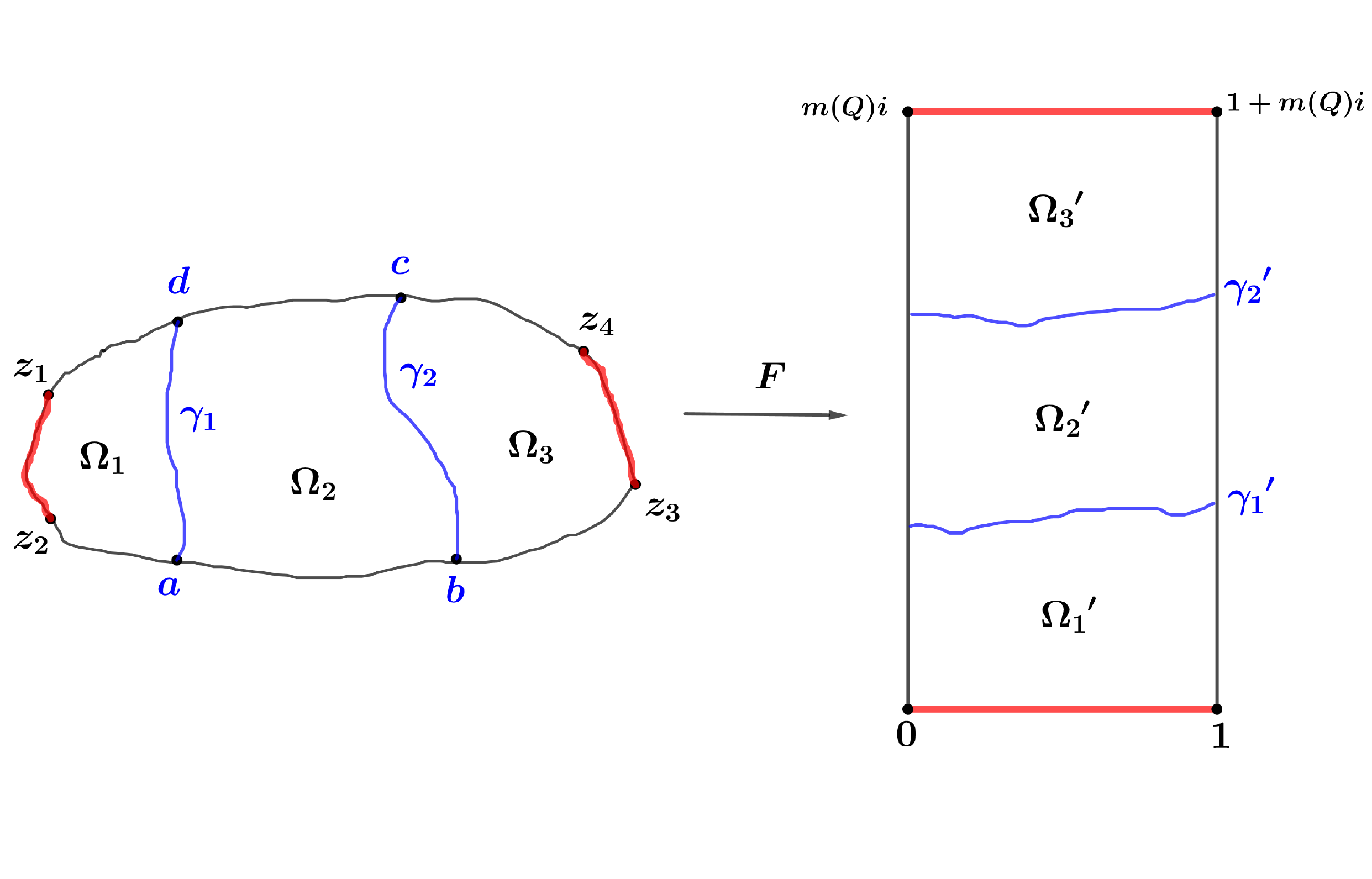}
		\caption{The subdivision of $\Omega$ into $\Omega_1,\Omega_2,\Omega_3$ and the conformal map $F:Q \to Q'$.}
		\label{dd}
	\end{center}
\end{figure}

Let $\Omega $ be a Jordan domain in $\mathbb{C}$ and consinder a system consisting of $\Omega $ and four distinct points $z_1,z_2,z_3,z_4$ in counterclockwise order on its boundary $\partial \Omega$. Such a system is said to be a quadrilateral $Q$ and is denoted by
\[Q: = \left\{ {\Omega ;{z_1},{z_2},{z_3},{z_4}} \right\}.\] 
The conformal module $m\left( Q \right)$ of $Q$ is the unique number for which $Q$ is conformally equivalent to the rectangular quadrilateral
\[Q': = \left\{ {R_{m\left( Q \right)} ;0,1,1+m\left( Q \right)i,m\left( Q \right)i} \right\},\] 
where $R_{m\left( Q \right)}=\left\{ {x + yi:0 < x < 1,0 < y < m\left( Q \right)} \right\}$ (see Fig. \ref{dd}). Note that $m\left( Q \right)$ is conformally invariant and it is equal to the extremal distance between the boundary arcs $\left( {z_1},{z_2} \right)$ and $\left( {z_3},{z_4} \right)$ of $\Omega$. So, $\Omega$ and $Q: = \left\{ {\Omega ;{z_1},{z_2},{z_3},{z_4}} \right\}$ will denote respectively the original domain and the corresponding quadrilateral. Moreover, ${\Omega _1},{\Omega _2}, \ldots ,$ and ${Q _1},{Q _2}, \ldots ,$ will denote the principle subdomains and corresponding component quadrilaterals of the decomposition under consideration. Now consider the situation of Fig. \ref{dd}, where the decomposition of $Q: = \left\{ {\Omega ;{z_1},{z_2},{z_3},{z_4}} \right\}$ is defined by two non-intersecting arcs $\gamma_1,\gamma_2$ that join respectively two distinct points $a$ and $b$ on the boundary arc $\left( {{z_2},{z_3}} \right)$ to two points $d$ and $c$ on the boundary arc $\left( {{z_4},{z_1}} \right)$. These two arcs subdivide $\Omega$ into three non-intersecting subdomains denoted by ${\Omega _1},{\Omega _2}$ and ${\Omega _3}$. In addition, the arc $\gamma_1$ subdivides $\Omega$ into $\Omega_1$ and another subdomain denoted by ${\Omega _{2,3}}$, i.e. we take
\[{\overline \Omega  _{2,3}} = {\overline \Omega  _2} \cup {\overline \Omega  _3}.\]
Similarly, we say that $\gamma_2$ subdivides $\Omega$ into $\Omega_{1,2}$ and $\Omega  _3$, i.e. we take
\[{\overline \Omega  _{1,2}} = {\overline \Omega  _1} \cup {\overline \Omega  _2}.\]
Finally, we use the notations $Q_1,Q_2,Q_3,Q_{1,2}$ and $Q_{2,3}$ to denote, respectively, the quadrilaterals corresponding to the subdomains $\Omega_1,\Omega_2,\Omega_3,\Omega_{1,2}$ and $\Omega_{2,3}$, i.e.
\[Q_1: = \left\{ {\Omega_1 ;{z_1},{z_2},a,d} \right\},\;Q_2: = \left\{ {\Omega_2 ;d,a,b,c} \right\},\;Q_3: = \left\{ {\Omega_3 ;c,b,z_3,z_4} \right\}\]
and
\[Q_{1,2}: = \left\{ {\Omega_{1,2} ;{z_1},{z_2},b,c} \right\},\;Q_{2,3}: = \left\{ {\Omega_{2,3} ;d,a,z_3,z_4} \right\}.\]

The following theorem was proved by Papamichael and Stylianopoulos in \cite[p. 221-222]{Pap}; see also \cite[p. 454-455]{Ga}.

\begin{theorem}\label{p1} Consider the decomposition and the notations illustrated in Fig. \ref{dd}. With the terminology  defined above, we have
	\[\left| {m\left( Q \right) - \left( {m\left( {{Q_{1,2}}} \right) + m\left( {{Q_{2,3}}} \right) - m\left( {{Q_{2}}} \right)} \right)} \right| \le 8.82{e^{ - \pi m\left( {{Q_2}} \right)}},\]
	provided that $m\left( {{Q_2}} \right) \ge 3$.
\end{theorem}

\begin{remark}\label{re} Papamichael and Stylianopoulos proved Theorem \ref{p1} in case $\Omega $ is a Jordan domain. However, it follows from the proof that the theorem is still valid if $\Omega $ is a simply connected domain and its boundary sets $\left( {{z_1},{z_2}} \right)$ and $\left( {{z_3},{z_4}} \right)$ are arcs of prime ends.
\end{remark}

\subsection{Beurling's estimates for harmonic measure}

A basic property of extremal distance is its connection to harmonic measure as the following theorems, due mainly to Beurling, state (see \cite[p. 280]{Bet}, \cite[p. 369-372]{Be}, \cite[p. 143-146]{Gar} and \cite[p. 100]{koo}).

\begin{figure}[H]
	\begin{center}
		\includegraphics[scale=0.27]{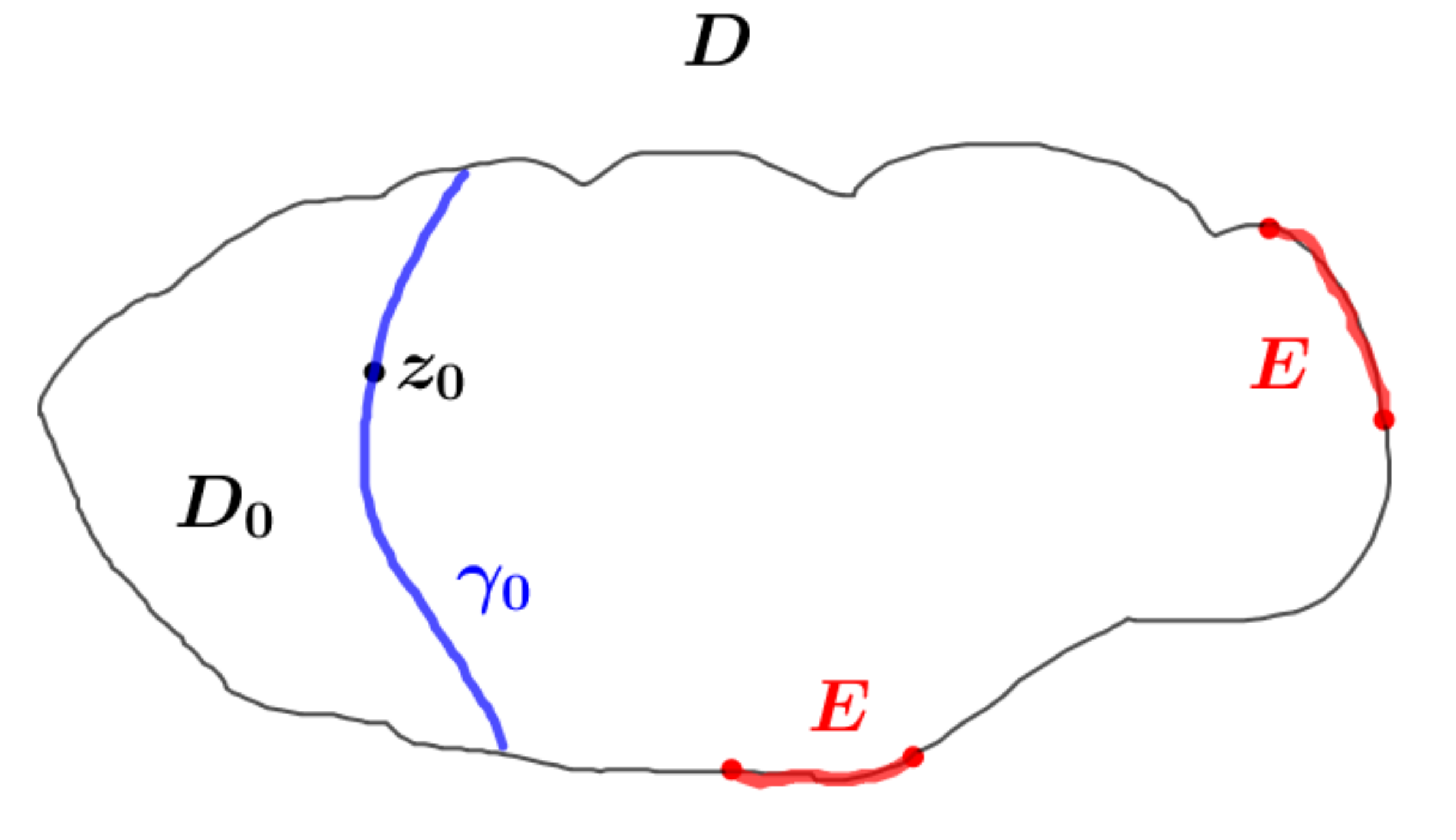}
		\caption{The simply connected domains $D$ and $D_0$.}
		\label{eikona}
	\end{center}
\end{figure}

\begin{theorem}\label{beu} Let $D$ be a simply connected domain in $\mathbb{C}$ and $E$ consist of a finite number of arcs lying on $\partial D$. Fix $z_0 \in D$ and choose a curve $\gamma_0$ that contains $z_0$, lies in $D$ and joins two points of $\partial D$ so that $\gamma_0$ bounds with $\partial D$ a domain $D_0$ and $z_0$ can be joined to $E$ inside $D\backslash {D_0}$ (see Fig. \ref{eikona}). If ${\lambda _{D\backslash {D_0}}}\left( {{\gamma _0},E} \right) > 2$, then
	\[{\omega _D}\left( {{z_0},E} \right) \le 3\pi {e^{ - \pi {\lambda _{D\backslash {D_0}}}\left( {{\gamma _0},E} \right)}}.\]
\end{theorem}

\begin{theorem}\label{ge} Let $D$ be a simply connected domain in $\mathbb{C}$ and $E$ be an arc (of prime ends) on $\partial D$. Fix $z_0 \in D$ and map $D$ onto $\mathbb{D}$ by the conformal map $f$ so that ${f\left( {{z_0}} \right) = 0}$ and $f\left( E \right) = \left\{ {{e^{i\theta }}:\theta  \in \left[ { - t,t} \right]} \right\}$ for some $t \in \left[ {0,\pi } \right]$. If ${\gamma _E}: = {f^{ - 1}}\left( {\left[ { - 1,0} \right]} \right)$, then there exists an absolute positive constant $C$ such that
	\[{\omega _D}\left( {{z_0},E} \right) \ge C{e^{ - \pi {\lambda _D}\left( {{\gamma _E},E} \right)}}.\]
\end{theorem} 
	
\section{Auxilary lemmas}\label{section3}

\begin{figure} [H]
	\begin{center}
		\includegraphics[scale=0.35]{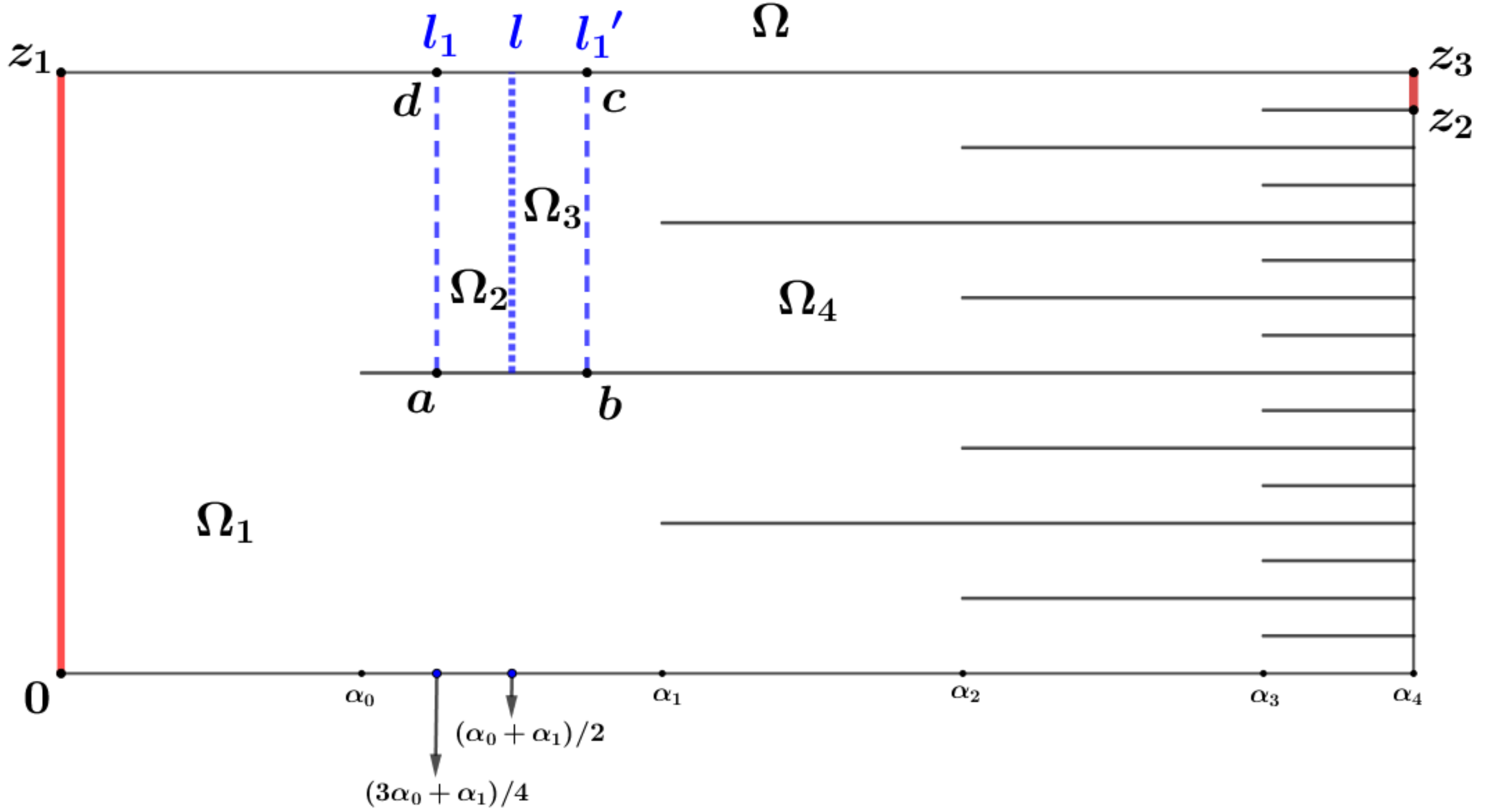}
		\caption{The decomposition of $\Omega$  into $\Omega_1,\Omega_2,\Omega_3,\Omega_4$.}
		\label{l11}
	\end{center}
\end{figure}

Let $\Omega$ be a simply connected domain of the form illustrated in Fig. \ref{l11}. Note that the positive numbers ${\alpha _0},{\alpha _1},{\alpha _2}, \ldots ,$ are the real parts of the tips of the horizontal boundary segments of $\Omega$. We consider the straight line crosscuts $l_1,l$ of Fig. \ref{l11} so that $l$ lies on the vertical line passing through the midpoint of $\left[ {{\alpha _0},{\alpha _1}} \right]$ and $l_1$ lies on the vertical line passing through the midpoint of $\left[ {{\alpha _0},\frac{{{\alpha _0} + {\alpha _1}}}{2}} \right]$. We decompose $\Omega$ by means of the straight line crosscuts $l_1,l,{l_1}'$ into four subdomains $\Omega_1,\Omega_2,\Omega_3,\Omega_4$ so that $\Omega_3$ is the reflection  of $\Omega_2$ in $l$.

\begin{lemma}\label{lem1} With the notation above, let $Q = \left\{ {\Omega ;{z_1},0,{z_2},{z_3}} \right\}$. According to the terminology introduced in Section \ref{section2}, we have for the decomposition defined by $l$,
	\[0 \le m\left( Q \right) - \left( {m\left( {{Q_{1,2}}} \right) + m\left( {{Q_{3,4}}} \right)} \right) \le 26.46{e^{ -\pi m\left( {{Q_2}} \right)}},\]
	provided that $m\left( {{Q_2}} \right) \ge 3$.
\end{lemma}

\proof Since $\left[ {{z_1},0} \right]$ and $\left[ {{z_2},{z_3}} \right]$ are arcs of prime ends, there exists a conformal map $F$ of $Q$ onto $F\left( Q \right) = \left\{ {{R_{m\left( Q \right)}};0,1,1 + m\left( Q \right)i,m\left( Q \right)i} \right\}$, where ${R_{m\left( Q \right)}} = \left\{ {x + yi:0 < x < 1,0 < y < m\left( Q \right)} \right\}$. By symmetry we have that ${m\left( {{Q_{2,3}}} \right) = 2m\left( {{Q_2}} \right) \ge 6}$. So, applying Theorem \ref{p1} we get
\[\left| {m\left( Q \right) - \left( {m\left( {{Q_{1,2,3}}} \right) + m\left( {{Q_{2,3,4}}} \right) - m\left( {{Q_{2,3}}} \right)} \right)} \right| \le 8.82{e^{ - \pi m\left( {{Q_{2,3}}} \right)}}\]
or equivalently
\begin{equation}\label{sx1}
\left| {m\left( Q \right) - \left( {m\left( {{Q_{1,2,3}}} \right) + m\left( {{Q_{2,3,4}}} \right) - 2m\left( {{Q_2}} \right)} \right)} \right| \le 8.82{e^{ - 2\pi m\left( {{Q_2}} \right)}}.
\end{equation}
Now consider the quadrilateral $Q_{2,3,4}$. Since $m\left( {{Q_3}} \right)=m\left( {{Q_2}} \right) \ge 3$, by applying Theorem \ref{p1}, we deduce that
\[\left| {m\left( Q_{2,3,4} \right) - \left( {m\left( {{Q_{2,3}}} \right) + m\left( {{Q_{3,4}}} \right) - m\left( {{Q_{3}}} \right)} \right)} \right| \le 8.82{e^{ - \pi m\left( {{Q_{3}}} \right)}}\]
or equivalently
\begin{equation}\label{sx2}
\left| {m\left( Q_{2,3,4} \right) - \left( {m\left( {{Q_{2}}} \right) + m\left( {{Q_{3,4}}} \right)} \right)} \right| \le 8.82{e^{ - \pi m\left( {{Q_{2}}} \right)}}.
\end{equation}
Similarly, consider the quadrilateral $Q_{1,2,3}$. Since $m\left( {{Q_2}} \right) \ge 3$, by applying Theorem \ref{p1}, we deduce that
\[\left| {m\left( Q_{1,2,3} \right) - \left( {m\left( {{Q_{1,2}}} \right) + m\left( {{Q_{2,3}}} \right) - m\left( {{Q_{2}}} \right)} \right)} \right| \le 8.82{e^{ - \pi m\left( {{Q_{2}}} \right)}}\]
or equivalently
\begin{equation}\label{sx3}
\left| {m\left( Q_{1,2,3} \right) - \left( {m\left( {{Q_{2}}} \right) + m\left( {{Q_{1,2}}} \right)} \right)} \right| \le 8.82{e^{ - \pi m\left( {{Q_{2}}} \right)}}.
\end{equation}
By relations (\ref{sx1}), (\ref{sx2}), (\ref{sx3}) and the serial rule
\[{m\left( Q \right) \ge m\left( {{Q_{1,2}}} \right) + m\left( {{Q_{3,4}}} \right)},\]
we finally get
\[0 \le m\left( Q \right) - \left( {m\left( {{Q_{1,2}}} \right) + m\left( {{Q_{3,4}}} \right)} \right) \le 26.46{e^{ -\pi m\left( {{Q_2}} \right)}}.\] 
\qed

In the following lemma we use the notation $D\left( {0,\alpha} \right)$ to denote the disk with center at 0 and radius $\alpha$.

\begin{figure}[H]
	\begin{center}
		\includegraphics[scale=0.35]{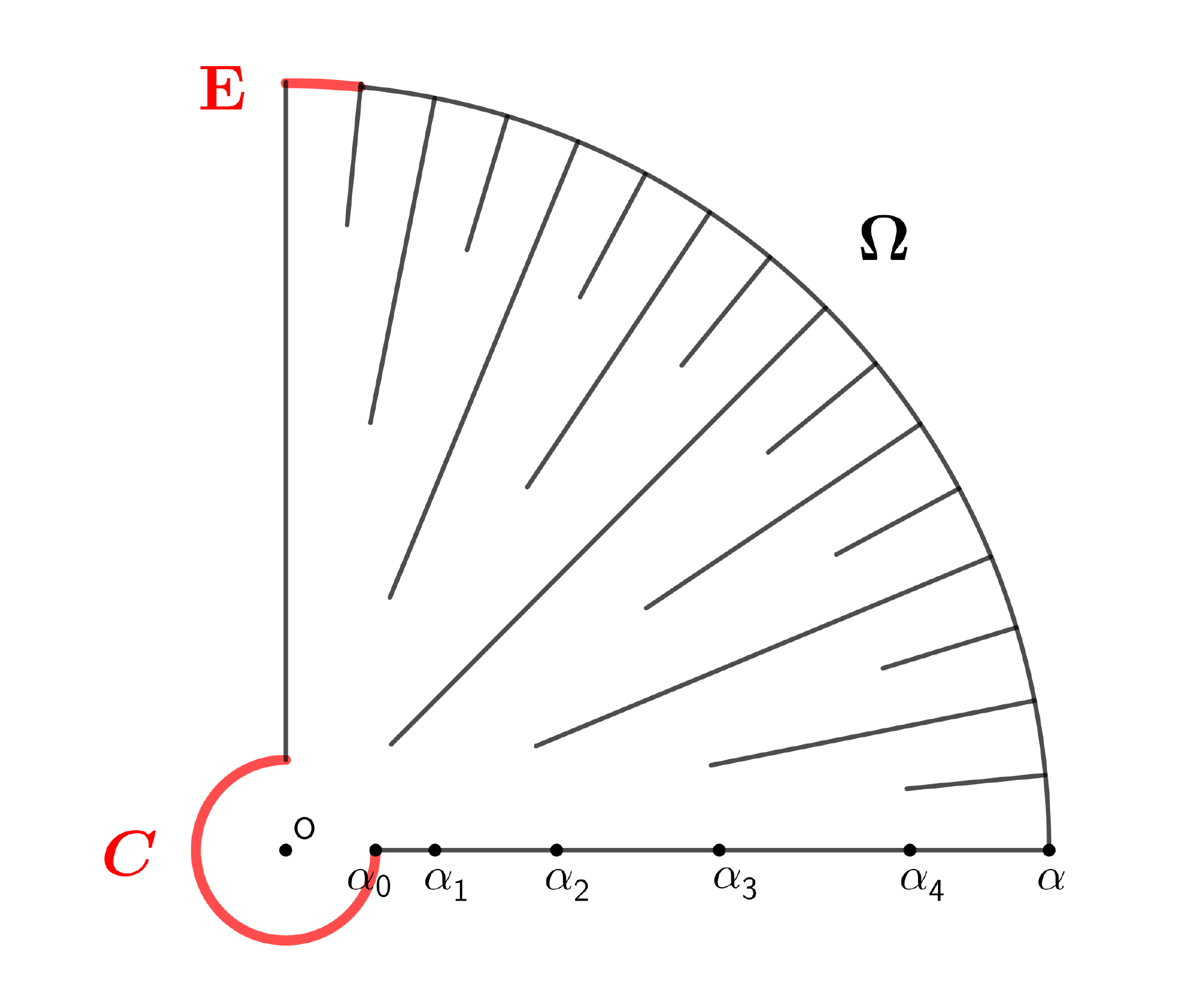}
		\caption{The simply connected domain $\Omega$ and the arc $E$.}
		\label{l21}
	\end{center}
\end{figure}

\begin{lemma}\label{lem2} Let $\Omega$ be a simply connected domain of the form illustrated in Fig. \ref{l21} and $E$ be an arc of prime ends on $\partial \Omega  \cap \partial D\left( {0,\alpha} \right)$. If $f$ is the conformal map of $\Omega$ onto $\mathbb{D}$ such that $f\left( 0 \right) = 0$ and $f\left( E \right) = \left\{ {{e^{i\theta }}:\theta  \in \left[ { - t,t} \right]} \right\}$ for some ${t \in \left[ {0,\pi } \right]}$, then 
	\[\gamma  \subset \overline D \left( {0,\alpha_0 } \right),\]
	where $\gamma : = {f^{ - 1}}\left( {\left[ { - 1,0} \right]} \right)$.
\end{lemma}

\begin{figure}[H]
	\begin{minipage}{0.4\textwidth}
		\begin{center}
			\includegraphics[scale=0.4]{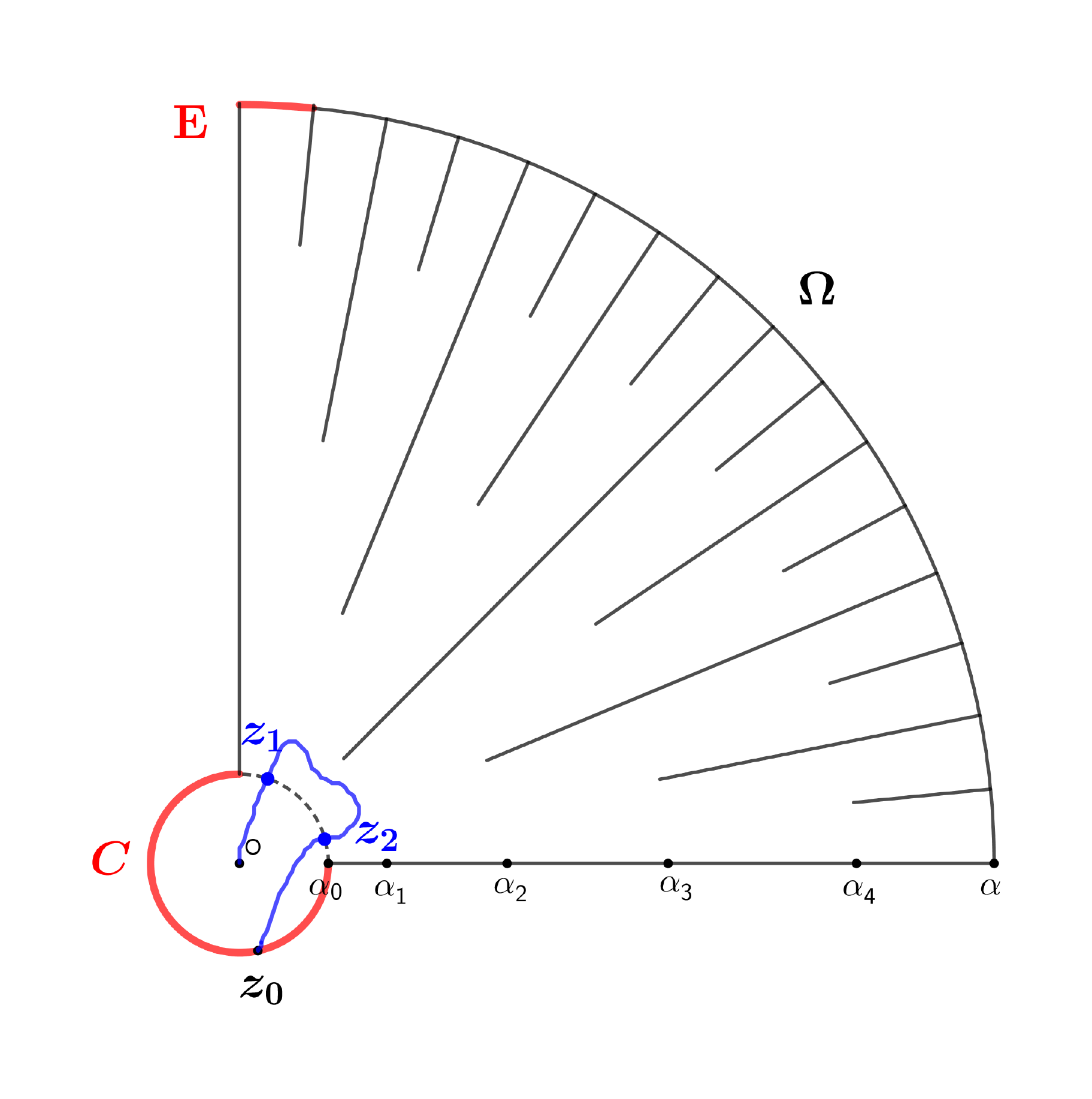}
			\caption{The curve $\gamma_0$ and the points $z_1,z_2$ in $\Omega$.}
			\label{l22}
		\end{center}
	\end{minipage}\hfill
	\begin{minipage}{0.6\textwidth}
		\begin{center}
			\includegraphics[scale=0.25]{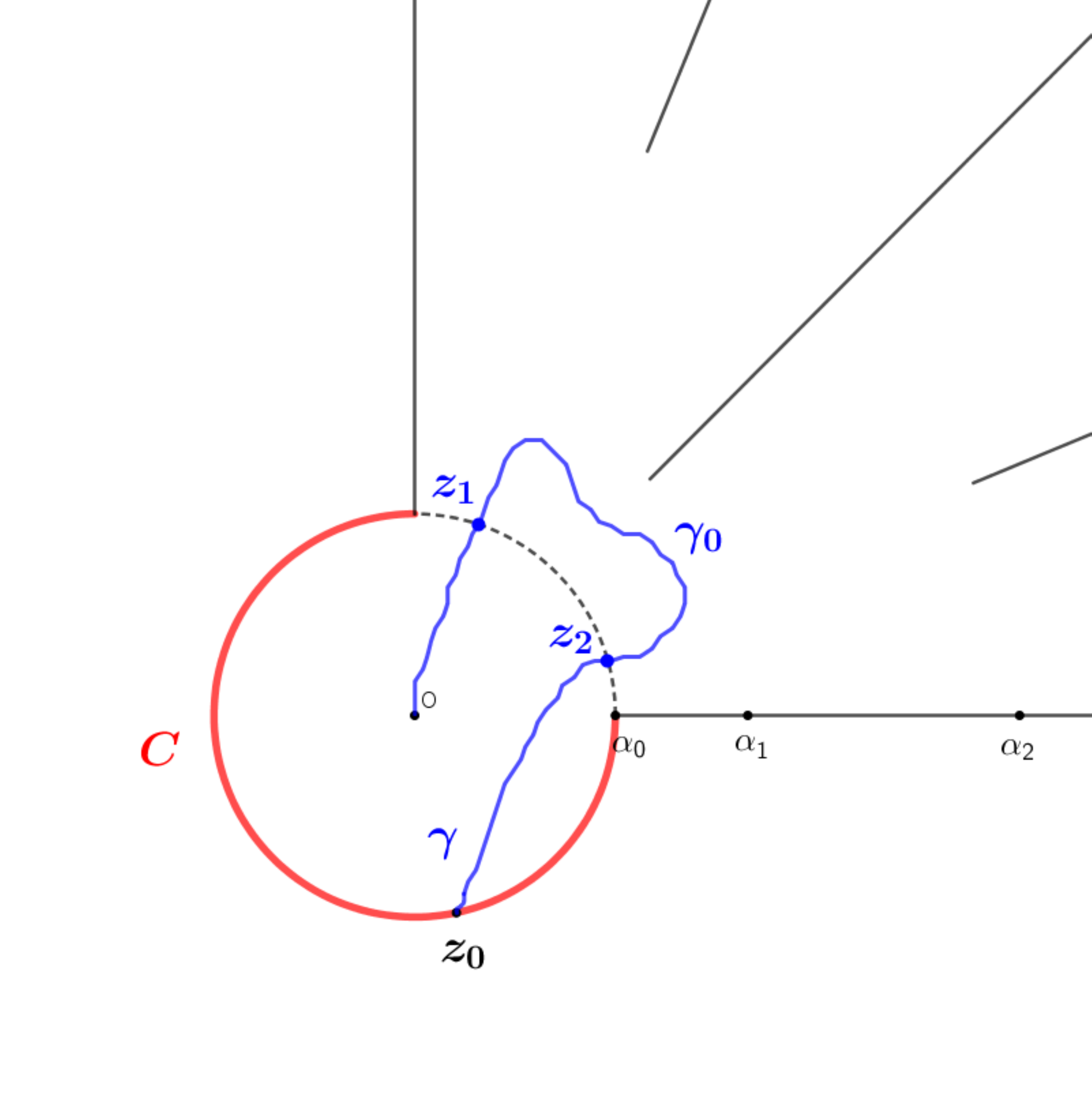}
			\caption{The curve $\gamma_0$ and the points $z_1,z_2$ in magnification.}
			\label{l23}
		\end{center}
	\end{minipage}
\end{figure}

\proof Set $C = \left\{ {\alpha_0 {e^{i\theta }}:\theta  \in \left[ {\frac{\pi }{2},2\pi } \right]} \right\}$ and ${z_0} = {f^{ - 1}}\left( { - 1} \right)$. Since $D\left( {0,\alpha_0 } \right) \subset \Omega$, by Corollary 4.3.9 \cite[p. 102]{Ra} and conformal invariance of harmonic measure, we have
\[{\omega _\mathbb{D}}\left( {0,f\left( C \right)} \right) = {\omega _\Omega }\left( {0,C} \right) \ge {\omega _{D\left( {0,\alpha_0 } \right)}}\left( {0,C} \right) = \frac{3}{4}.\]
This, in conjunction with the fact that $f\left( C \right) \cap f\left( E \right) = \emptyset $ and $f\left( C \right)$ is a connected arc of $\partial \mathbb{D}$, implies that
\[\left\{ {{e^{i\theta }}:\theta  \in \left[ {\frac{\pi }{2},\frac{{3\pi }}{2}} \right]} \right\} \subset f\left( C \right).\]
So, ${z_0} = {f^{ - 1}}\left( { - 1} \right) \in C$. Now suppose that $\gamma  \not\subset \overline D \left( {0,\alpha_0 } \right)$. Then $\gamma$ contains a curve $\gamma_0$ lying in $\Omega \backslash D\left( {0,\alpha_0 } \right)$ with endpoints $z_1,z_2 \in \partial {D \left( {0,\alpha_0 } \right)}$ (see Fig. \ref{l22}, \ref{l23}).

Since $\gamma$ is the hyperbolic geodesic joining $0$ to ${z_0}$ in $\Omega$, $\gamma_0$ is the hyperbolic geodesic joining $z_1$ to ${z_2}$ in $\Omega$. Notice that $\Omega$ is a hyperbolic region in $\mathbb{C}$ such that $\Omega  \cap \partial D\left( {0,\alpha_0 } \right) \ne \emptyset $ and $\Omega \backslash D\left( {0,\alpha_0 } \right) \subset \Omega ^ *$, where $\Omega ^ *$ is the reflection of $\Omega$ in the circle $\partial D\left( {0,\alpha_0 } \right)$. So, applying Theorem \ref{rp} we get
\[{\lambda _{{\Omega ^ * }}}\left( z \right) < {\lambda _\Omega }\left( z \right),\;z \in {\gamma _0}\]
and thus
\[\int_{{\gamma _0}^ * } {{\lambda _\Omega }\left( {{z^ * }} \right)\left| {d{z^ * }} \right|}  < \int_{{\gamma _0}} {{\lambda _\Omega }\left( z \right)\left| {dz} \right|}, \]
where ${\gamma _0}^ * $ is the reflection of $\gamma _0$ in $\partial D\left( {0,\alpha_0 } \right)$. But this leads to contradiction because $\gamma_0$ is the hyperbolic geodesic joining $z_1$ to ${z_2}$ in $\Omega$. So,
$\gamma  \subset \overline D \left( {0,\alpha_0 } \right)$. Note that the same result could come from J\o rgensen' s theorem that closed disks in $\Omega$ are strictly convex in the hyperbolic geometry of $\Omega$ (see \cite{Jo}).

\qed 

\begin{figure}[H]
	\begin{center}
		\includegraphics[scale=0.48]{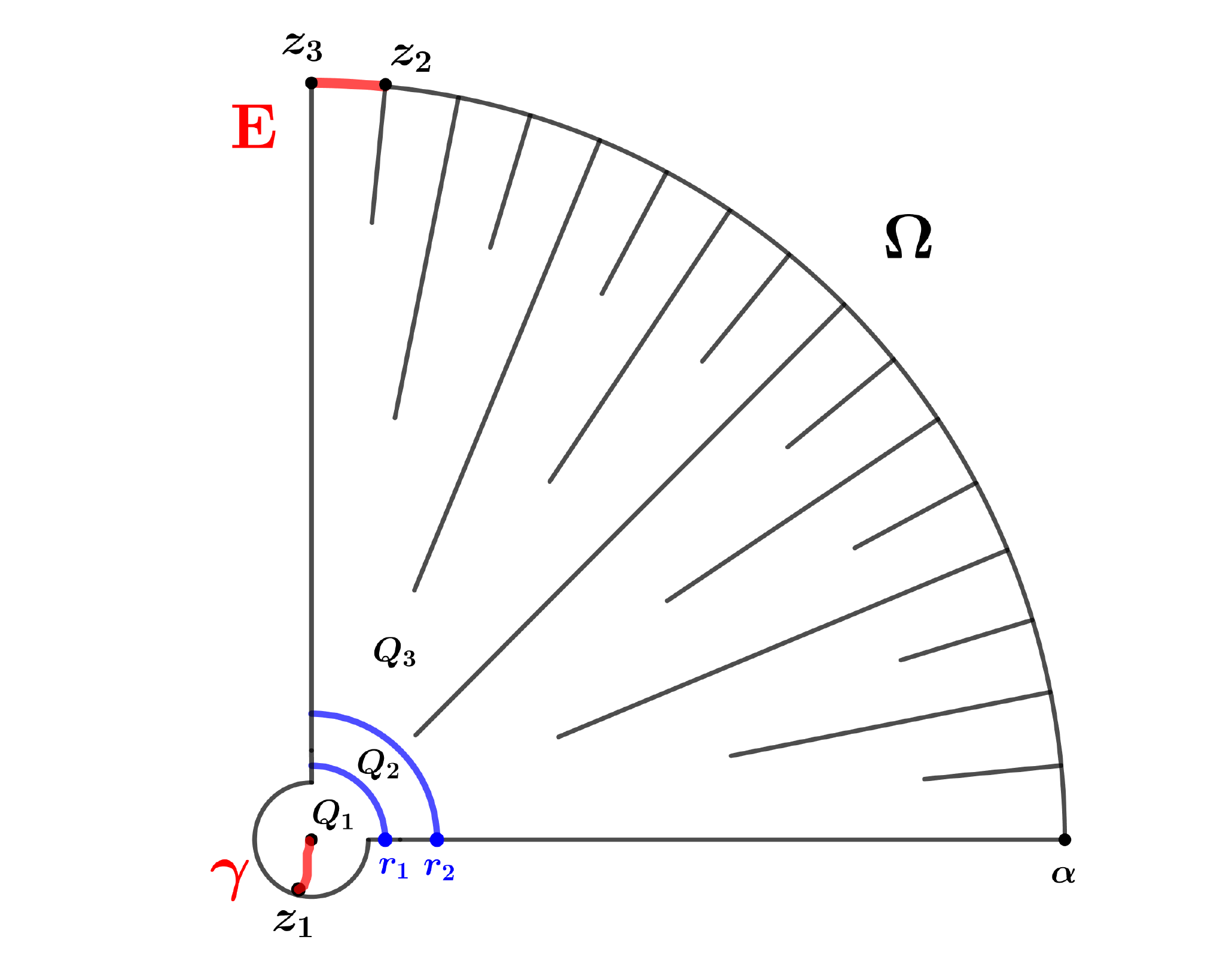}
		\caption{The decomposition of $Q$ into $Q_1,Q_2,Q_3$.}
		\label{l31}
	\end{center}
\end{figure}

\begin{lemma}\label{lem3} Let $\Omega,\gamma,E$ be as in Lemma \ref{lem2} and ${z_1},{z_2},{z_3}$ be the points illustrated in Fig. \ref{l31}. Take $r_1,r_2$ so that $\alpha_0  < {r_1} < {r_2} < {\alpha _1}$ and $\log \frac{{{r_2}}}{{{r_1}}} \ge \frac{{3\pi }}{2}$. Decomposing $Q = \left\{ {{\Omega \backslash \gamma } ;{z_1},{z_1},{z_2},{z_3}} \right\}$ as in Fig. \ref{l31}, with the terminology introduced in Section \ref{section2}, we have
	\[\left| {m\left( Q \right) - \left( {m\left( {{Q_{1,2}}} \right) + m\left( {{Q_{2,3}}} \right) - m\left( {{Q_2}} \right)} \right)} \right| \le 8.82{e^{ - \pi m\left( {{Q_2}} \right)}}.\]
In the notation $Q = \left\{ {{\Omega \backslash \gamma } ;{z_1},{z_1},{z_2},{z_3}} \right\}$, by the pair of points ${z_1},{z_1}$, we mean the two different prime ends supported at the point ${z_1}$. 
\end{lemma}

\proof First, applying the conformal map $f\left( z \right) = \LOG z$ (principal branch of the logarithm) on the quadrilateral ${Q_2} = \left\{ {{\Omega _2};{r_1}i,{r_1},{r_2},{r_2}i} \right\}$ we take the rectangular quadrilateral
\[f\left( {{Q_2}} \right) = \left\{ {f\left( {{\Omega _2}} \right);\log{r_1} + \frac{\pi }{2}i,\log {r_1},\log{r_2},\log{r_2} + \frac{\pi }{2}i} \right\},\]
where $f\left( {{\Omega _2}} \right) = \left\{ {x + yi:\log {r_1} < x < \log {r_2},0 < y < {\pi  \mathord{\left/
			{\vphantom {\pi  2}} \right.
			\kern-\nulldelimiterspace} 2}} \right\}$. Because of the conformal invariance of modules and our assumption about $r_1,r_2$,
\[m\left( {{Q_2}} \right) = m\left( f\left( {{Q_2}} \right) \right) = \frac{{\log \left( {{{{r_2}} \mathord{\left/
					{\vphantom {{{r_2}} {{r_1}}}} \right.
					\kern-\nulldelimiterspace} {{r_1}}}} \right)}}{{{\pi  \mathord{\left/
				{\vphantom {\pi  2}} \right.
				\kern-\nulldelimiterspace} 2}}} \ge 3.\]
Since the boundary sets $\left( {{z_1},0,{z_1}} \right)$ and $E$ are arcs of prime ends, there exists a conformal map $F$ of $Q$ onto 
\[Q' = \left\{ {R_{m\left( Q \right)};0,1,1 + m\left( Q \right)i,m\left( Q \right)i} \right\},\]
where $R_{m\left( Q \right)} = \left\{ {x + yi:0 < x < 1,0 < y < m\left( Q \right)} \right\}$. Since $m\left( {{Q_2}} \right) \ge 3$, Theorem \ref{p1} implies that
\[\left| {m\left( Q \right) - \left( {m\left( {{Q_{1,2}}} \right) + m\left( {{Q_{2,3}}} \right) - m\left( {{Q_2}} \right)} \right)} \right| \le 8.82{e^{ - \pi m\left( {{Q_2}} \right)}}.\]
\qed

\section{The first example}\label{section4}

\begin{theorem}\label{main} There exists an unbounded simply connected domain $D$ with the following properties: Let $\psi $ be a conformal map of $\mathbb{D}$ onto $D$ with $\psi \left( 0 \right) = 0$. If ${F_\alpha } = \left\{ {z \in \mathbb{D}:\left| {\psi \left( z \right)} \right| = \alpha } \right\}$ for $\alpha >0$, then
\begin{enumerate}
\item the number of components of $\psi \left( {F_\alpha } \right)$ tends to infinity as $\alpha \to  + \infty$ and
\item $\forall K>0$ $\exists \alpha$ such that
\[{\omega _\mathbb{D}}\left( {0,{F_\alpha }} \right) \ge K{e^{ - {d_\mathbb{D}}\left( {0,{F_\alpha }} \right)}}.\]
\end{enumerate} 
\end{theorem}

\begin{figure}[H]
\begin{center}
\includegraphics[scale=0.33]{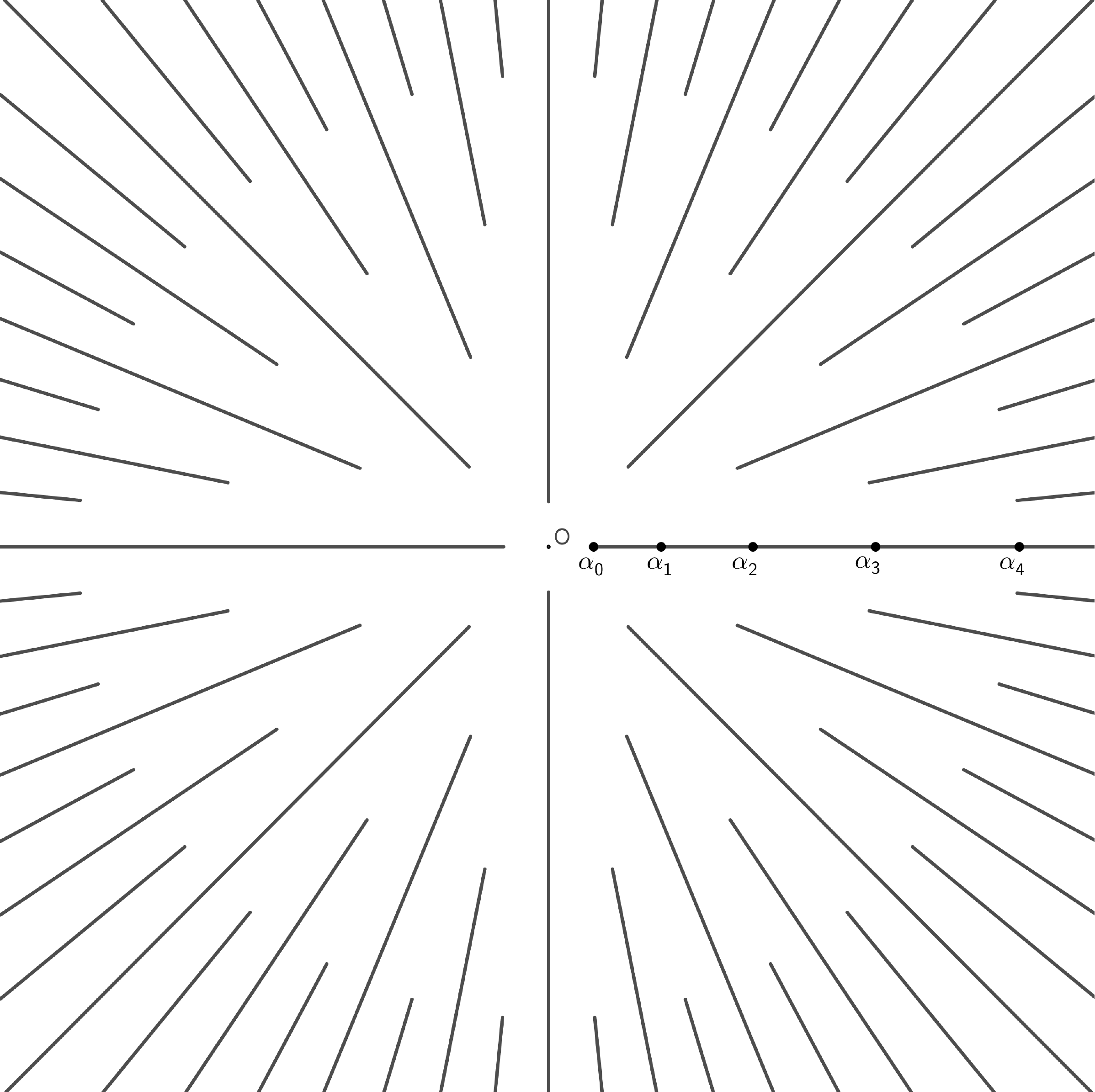}
\vspace*{0.3cm}
\caption{The simply connected domain $D$.}
\label{topos}
\end{center}
\end{figure}

\proof \textbf{Step 1:} If $\alpha_0=1,\;\alpha_1={e^{4\pi }},\;\alpha_2={e^{8\pi }},\ldots ,\alpha_n={e^{4n\pi }}, \ldots $, let $D$ be the simply connected domain of Fig. \ref{topos}, namely
\[D=\mathbb{C} \backslash \bigcup\limits_{k = 0}^3 {\left[ {{\alpha _0}{e^{i\frac{{k\pi }}{2}}}, + \infty } \right)} \backslash \bigcup\limits_{l = 1}^{ + \infty } {\bigcup\limits_{k = 0}^{{2^{l + 1}} - 1} {\left[ {{\alpha _l}{e^{i\frac{\pi }{{{2^l}}}\left( {\frac{1 }{2} + k} \right)}}, + \infty } \right)} },\]
with the notation $\left[ {r{e^{i\theta }}, + \infty } \right) = \left\{ {s{e^{i\theta }}:s \ge r} \right\}$.  

\begin{figure}[H]
\begin{minipage}{0.5\textwidth}
\begin{center}
\includegraphics[scale=0.5]{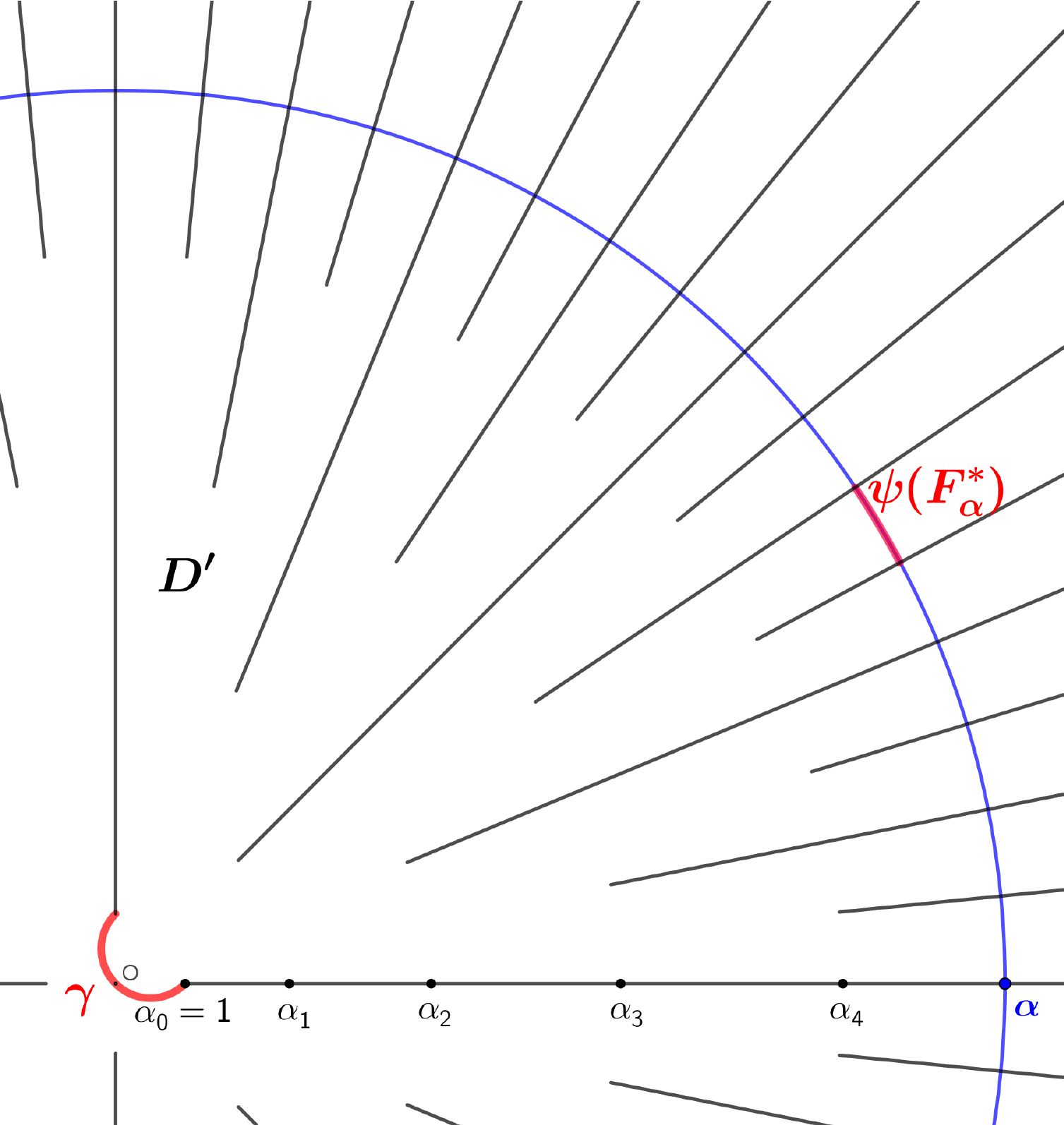}
\vspace*{0.3cm}
\caption{The simply connected domain $D'$ in case $\alpha  \in \left( {{\alpha _4},{\alpha _5}} \right)$.}
\label{d'}
\end{center}
\end{minipage}\hfill
\begin{minipage}{0.5\textwidth}
\begin{center}
\includegraphics[scale=0.5]{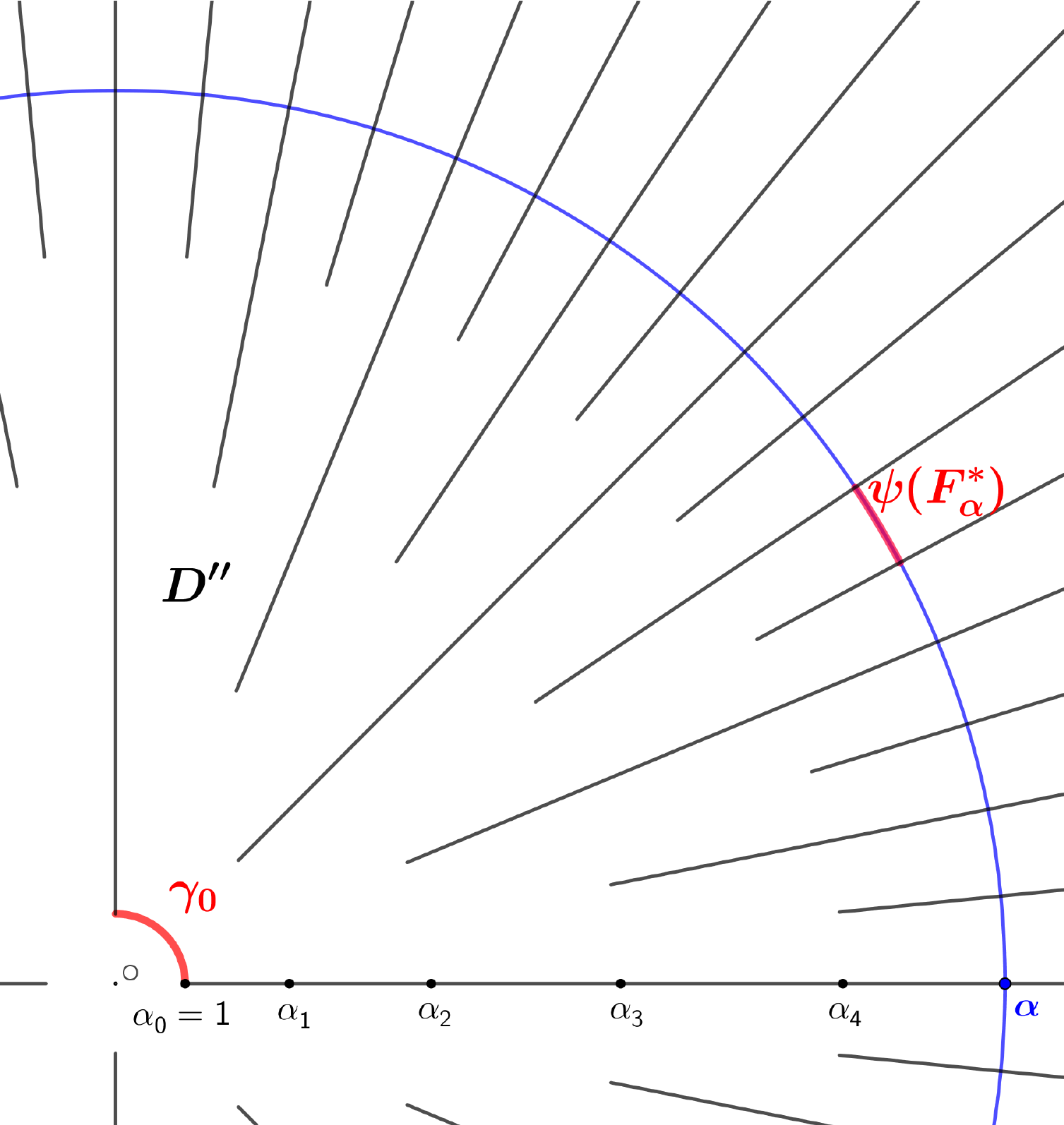}
\vspace*{0.3cm}
\caption{The simply connected domain $D''$ in case $\alpha  \in \left( {{\alpha _4},{\alpha _5}} \right)$.}
\label{d''}
\end{center}
\end{minipage}
\end{figure}

The Riemann Mapping Theorem implies that there exists a conformal map $\psi $ from $\mathbb{D}$ onto $D$ such that $\psi \left( 0 \right) = 0$. Let $N\left( \alpha  \right)$ be the number of components of $\psi \left( {{F_\alpha }} \right) = D \cap \partial D\left( {0,\alpha} \right)$, then we have
\[N\left( \alpha  \right) = \left\{ \begin{array}{l}
1,\;\;\;\;\;\;\;if\;\alpha  \in \left( {0,{\alpha _0}} \right) \\
{2^2},\;\;\;\;\;if\;\alpha  \in \left[ {{\alpha _0},{\alpha _1}} \right) \\
{2^3},\;\;\;\;\;if\;\alpha  \in \left[ {{\alpha _1},{\alpha _2}} \right)\\
 \vdots \\
{2^{n + 2}},\;if\;\alpha  \in \left[ {{\alpha _n},{\alpha _{n + 1}}} \right) \\
 \vdots 
\end{array} \right.\]
\textbf{Step 2:}  We fix a real number $\alpha$ such that $\alpha  > {\alpha _1}$ and $\alpha  \ne {\alpha _n}$ for every $n \in \mathbb{N}$. Then there exists a fixed number $n \in \mathbb{N}$ such that $\alpha  \in \left( {{\alpha _n},{\alpha _{n + 1}}} \right)$ and thus $N\left( \alpha  \right) = {2^{n + 2}}$. Since hyperbolic distance is conformally invariant we have
\begin{equation}\label{pr1}
{e^{ - {d_\mathbb{D}}\left( {0,{F_\alpha }} \right)}} = {e^{ - {d_D}\left( {0,\psi \left( {{F_\alpha }} \right)} \right)}} = {e^{ - {d_D}\left( {0,{\psi \left( {F_\alpha ^ * } \right)}} \right)}},
\end{equation}
where ${\psi \left( {F_\alpha ^ * } \right)}$ is a component of $\psi \left( {{F_\alpha }} \right)$ containing a point $z_0$ for which
\[{d_D}\left( {0,\psi \left( {{F_\alpha }} \right)} \right) = \inf \left\{ {{d_D}\left( {0,z} \right):z \in \psi \left( {{F_\alpha }} \right)} \right\} = {d_D}\left( {0,{z_0}} \right).\]
Due to the symmetry of $D$, we may assume without loss of generality that ${\psi \left( {F_\alpha ^ * } \right)}$  lies on the first quartile $P = \left\{ {z \in \mathbb{C}:\IM z > 0,\;\RE z > 0} \right\}$. By relation (\ref{1.1}) of Section \ref{section1}, we infer that
\[{e^{ - {d_\mathbb{D}}\left( {0,{F_\alpha ^ * } } \right)}} \le \frac{\pi }{2}{\omega _\mathbb{D}}\left( {0,{F_\alpha ^ * }} \right)\]
which in conjunction with the conformal invariance gives
\begin{equation}\label{pr2}
{e^{ - {d_D}\left( {0,{\psi \left( {F_\alpha ^ * } \right)}} \right)}}  \le  \frac{\pi }{2}{\omega _D}\left( {0,{\psi \left( {F_\alpha ^ * } \right)}} \right).
\end{equation}
Moreover, by Theorem \ref{beu} we deduce that
\begin{equation}\label{pr3}
{\omega _D}\left( {0,{\psi \left( {F_\alpha ^ * } \right)}} \right) \le 3\pi {e^{ - \pi {\lambda _{D'}}\left( {\gamma ,{\psi \left( {F_\alpha ^ * } \right)}} \right)}},
\end{equation}
where $\gamma$ is the arc of the circle passing through the points $0,\;1,\;i$ such that $\gamma$ connects $i$ to $1$ and $\gamma  \cap P = \emptyset$ and $D'$ is the subdomain of $D$ bounded by $\gamma$, $\left[ {1,\alpha } \right],\;\left[ {i,\alpha i} \right],\;P \cap \partial D\left( {0,\alpha} \right)$ and $\partial D \cap P \cap D\left( {0,\alpha} \right)$ (see Fig. \ref{d'}). If $\gamma_0=\partial \mathbb{D} \cap \overline P $ and $D''$ is the subdomain of $D'$ bounded by $\gamma_0$, $\left[ {1,\alpha } \right],\;\left[ {i,\alpha i} \right]\;P \cap \partial D\left( {0,\alpha} \right)$ and $\partial D \cap P \cap D\left( {0,\alpha} \right)$ (see Fig. \ref{d''}), then by Theorem \ref{el} we have
\begin{equation}\label{pr4}
{\lambda _{D'}}\left( {\gamma ,{\psi \left( {F_\alpha ^ * } \right)}} \right) \ge {\lambda _{D''}}\left( {{\gamma _0},{\psi \left( {F_\alpha ^ * } \right)}} \right).
\end{equation}
Combining the relations (\ref{pr1}), (\ref{pr2}), (\ref{pr3}) and (\ref{pr4}), we obtain
\begin{equation}\label{pr5}
{e^{ - {d_\mathbb{D}}\left( {0,{F_\alpha }} \right)}} \le \frac{{3{\pi ^2}}}{2}{e^{ - \pi {\lambda _{D''}}\left( {{\gamma _0},{\psi \left( {F_\alpha ^ * } \right)}} \right)}}.
\end{equation}
 
\begin{figure}[H]
\begin{center}
\includegraphics[scale=0.5]{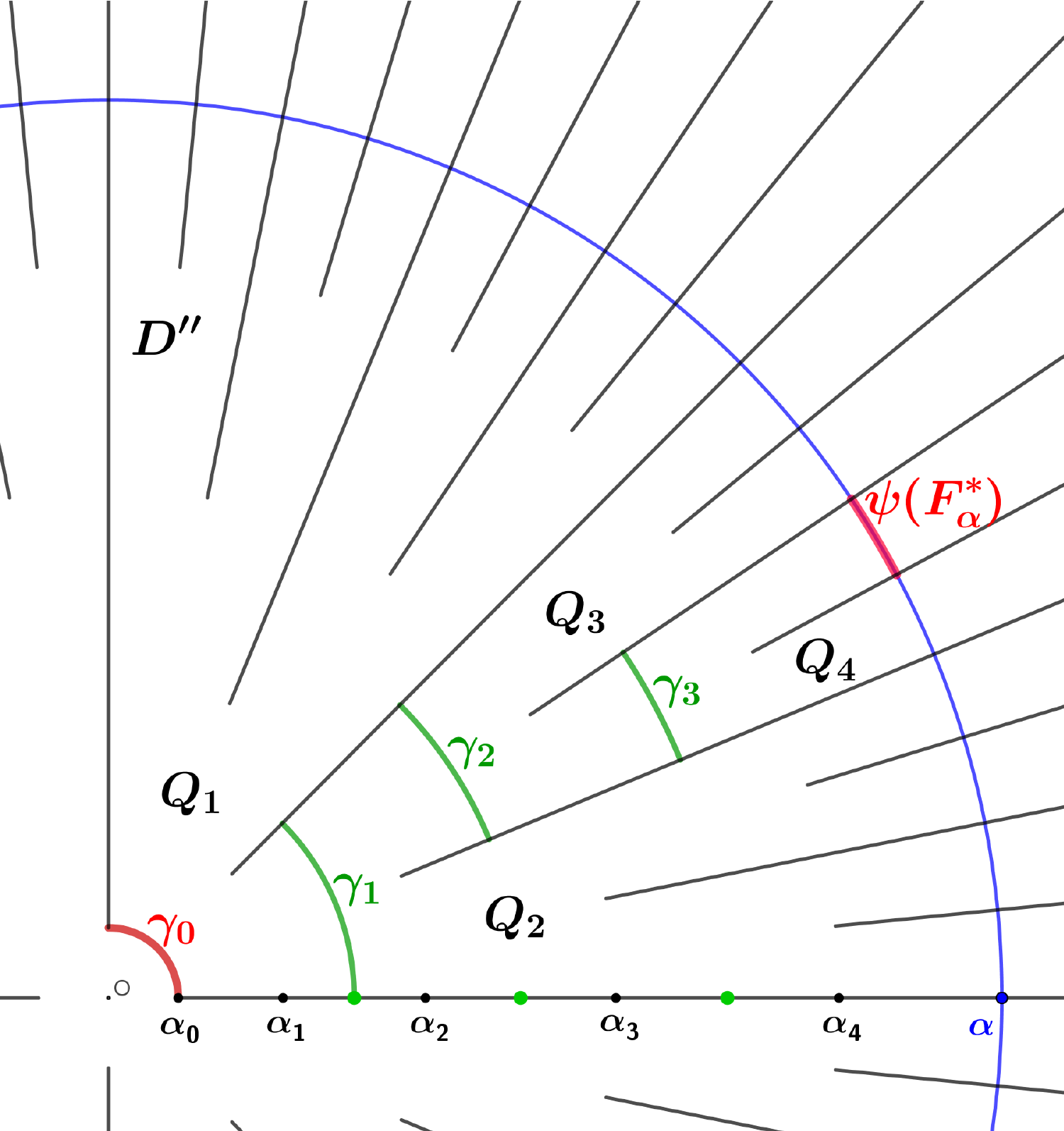}
\vspace*{0.3cm}
\caption{The crosscuts $\gamma _j$ of $D''$ in case $\alpha  \in \left( {{\alpha _4},{\alpha _5}} \right)$.}
\label{cross}
\end{center}
\end{figure}

Next we consider the crosscuts ${\gamma _1},{\gamma _2},{\gamma _3}, \ldots ,{\gamma _{n - 1}}$ of $D''$, where for $j = 1,2,3, \ldots ,n - 1$, $\gamma _j$ is an arc of the circle with center at $0$ and radius equal to the midpoint of ${\left[ {{\alpha _j},{\alpha _{j + 1}}} \right]}$ as illustrated in Fig. \ref{cross}.

Setting
\[m\left( {{Q_1}} \right): = {\lambda _{D''}}\left( {{\gamma _0},{\gamma _1}} \right),\;m\left( {{Q_2}} \right): = {\lambda _{D''}}\left( {{\gamma _1},{\gamma _2}} \right), \ldots,\;m\left( {{Q_n}} \right): = {\lambda _{D''}}\left( {{\gamma _{n - 1}},{\psi \left( {F_\alpha ^ * } \right)}} \right),\]
by the serial rule we deduce that
\begin{equation}\label{pr6}
{\lambda _{D''}}\left( {{\gamma _0},{\psi \left( {F_\alpha ^ * } \right)}} \right) \ge m\left( {{Q_1}} \right) + m\left( {{Q_2}} \right) + \ldots  + m\left( {{Q_n}} \right)
\end{equation}
and thus by (\ref{pr5}) and (\ref{pr6}),
\begin{equation}\label{pr7}
{e^{ - {d_\mathbb{D}}\left( {0,{F_\alpha }} \right)}} \le \frac{{3{\pi ^2}}}{2}{e^{ - \pi \left( {m\left( {{Q_1}} \right) + m\left( {{Q_2}} \right) +  \ldots  + m\left( {{Q_n}} \right)} \right)}}.
\end{equation}

\begin{figure}[H]
\begin{center}
\includegraphics[scale=0.52]{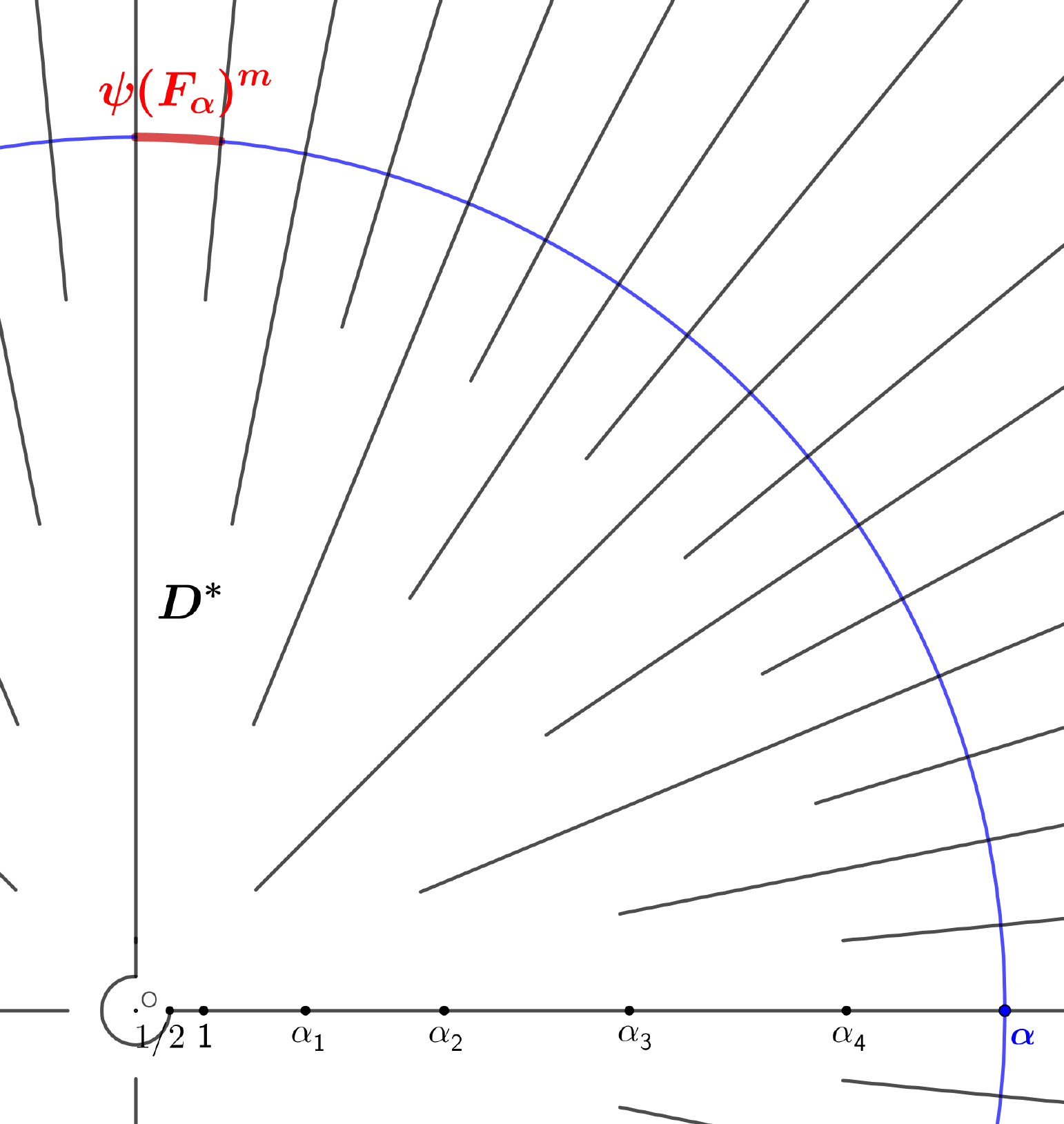}
\vspace*{0.3cm}
\caption{The simply connected domain $D^*$ in case $\alpha  \in \left( {{\alpha _4},{\alpha _5}} \right)$.}
\label{top}
\end{center}
\end{figure}

\textbf{Step 3:}  Since harmonic measure is conformally invariant and $N\left( \alpha  \right) = {2^{n + 2}}$, we have 
\[{\omega _\mathbb{D}}\left( {0,{F_\alpha }} \right) = {\omega _D}\left( {0,\psi \left( {{F_\alpha }} \right)} \right) = \sum\limits_{j = 1}^{N\left( \alpha  \right)} {{\omega _D}\left( {0,\psi {{\left( {{F_\alpha }} \right)}^j}} \right)}  \ge N\left( \alpha  \right){\omega _D}\left( {0,\psi {{\left( {{F_\alpha }} \right)}^m}} \right) = {2^{n + 2}}{\omega _D}\left( {0,\psi {{\left( {{F_\alpha }} \right)}^m}} \right),\]
where ${\psi {{\left( {{F_\alpha }} \right)}^j}}$, $j = 1,2,3, \ldots ,N(\alpha )$, are the components of ${\psi {{\left( {{F_\alpha }} \right)}}}$ and 
\[{\omega _D}\left( {0,\psi {{\left( {{F_\alpha }} \right)}^m}} \right) = \min \left\{ {{\omega _D}\left( {0,\psi {{\left( {{F_\alpha }} \right)}^j}} \right):j \in \left\{ {1,2,3, \ldots ,N(\alpha )} \right\}} \right\}.\] 
If $D^*$ is the subdomain of $D$ bounded by $\left[ {\frac{1}{2},\alpha } \right],\;\left[ {\frac{1}{2}i,\alpha i} \right],\;\partial D\left( {0,\frac{1}{2}} \right)\backslash P,\;\partial D\left( {0,\alpha } \right) \cap P$ and $\partial D \cap D\left( {0,\alpha } \right) \cap P$ as illustrated in Fig. \ref{top}, then applying Corollary 4.3.9 \cite[p. 102]{Ra} we obtain
\[{\omega _D}\left( {0,\psi {{\left( {{F_\alpha }} \right)}^m}} \right) \ge {\omega _{{D^*}}}\left( {0,\psi {{\left( {{F_\alpha }} \right)}^m}} \right)\]
and hence
\begin{equation}\label{pr8}
{\omega _\mathbb{D}}\left( {0,{F_\alpha }} \right) \ge {2^{n + 2}}{\omega _{{D^*}}}\left( {0,\psi {{\left( {{F_\alpha }} \right)}^m}} \right).
\end{equation}

\begin{figure}[H]
 \begin{minipage}{0.5\textwidth}
 \begin{center}
 \includegraphics[scale=0.27]{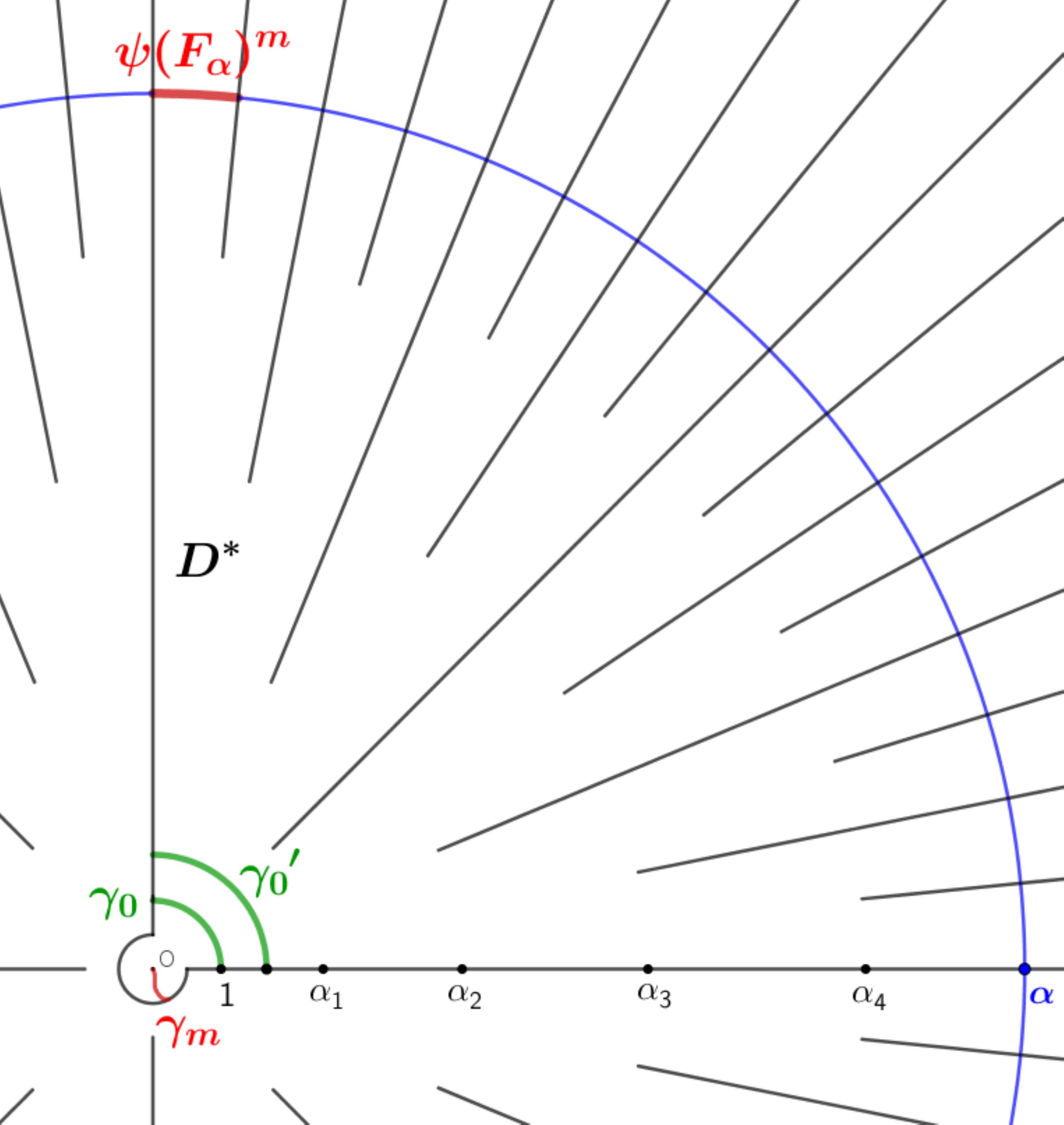}
 \vspace*{0.3cm}
 \caption{The crosscuts ${\gamma _0},{\gamma _0}'$ of $D^*$ and $\gamma _m$ in case $\alpha  \in \left( {{\alpha _4},{\alpha _5}} \right)$.}
 \label{top1}
 \end{center}
 \end{minipage}\hfill
 \begin{minipage}{0.5\textwidth}
 \begin{center}
 \includegraphics[scale=0.27]{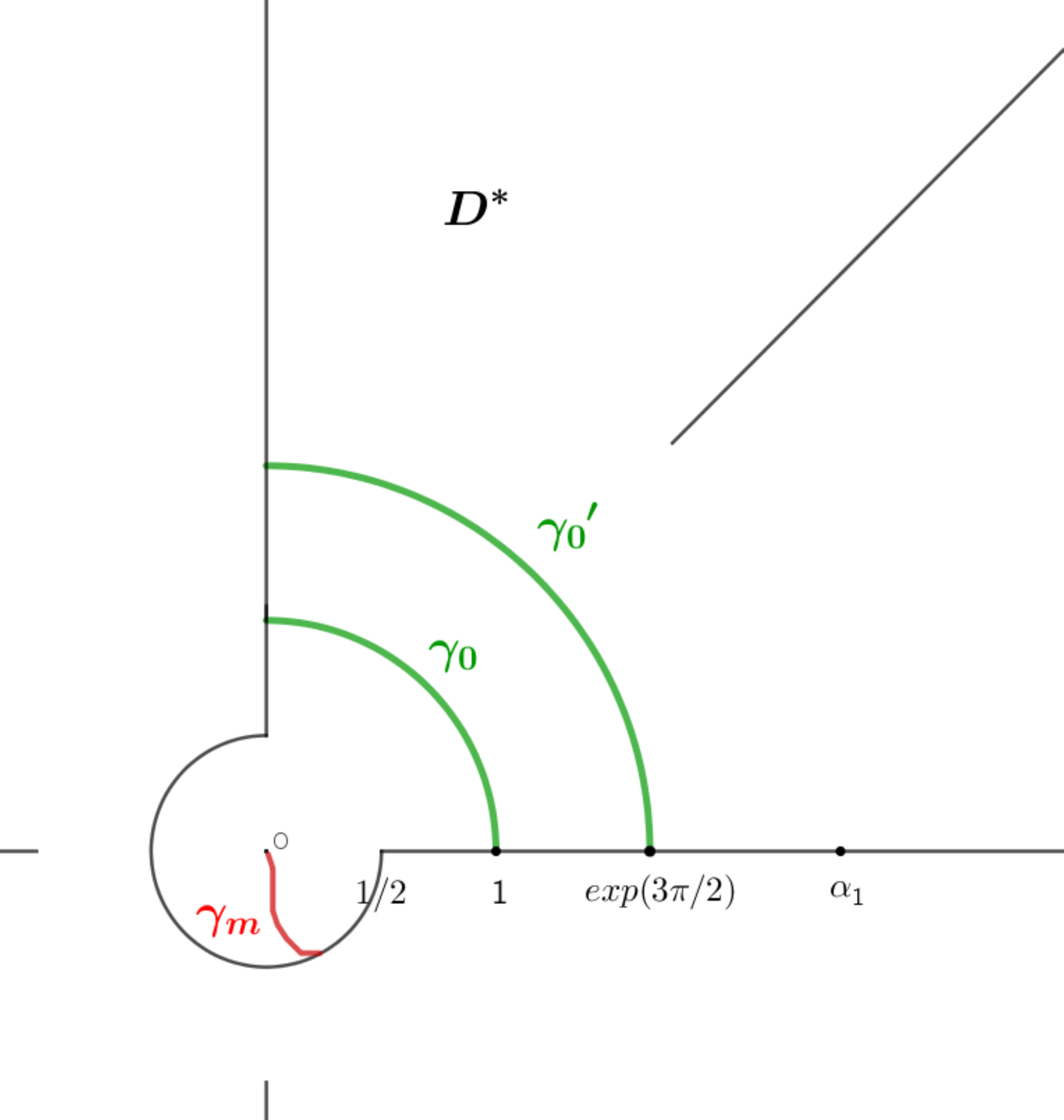}
 \vspace*{0.3cm}
 \caption{The crosscuts ${\gamma _0},{\gamma _0}'$ of $D^*$ and $\gamma _m$ in magnification.}
 \label{top2}
 \end{center}
 \end{minipage}
 \end{figure}

\textbf{Step 4:} If $h_m$ is the conformal map of $D^*$ onto $\mathbb{D}$ such that ${h_m}\left( 0 \right) = 0$ and ${h_m}\left( {\psi {{\left( {{F_\alpha }} \right)}^m}} \right) = \left\{ {{e^{i\theta }}:\theta  \in \left[ { - t,t} \right]} \right\}$ for some $t \in \left[ {0,\pi } \right]$ and ${\gamma _m}: = h_m^{ - 1}\left( {\left[ { - 1,0} \right]} \right)$, then by Theorem \ref{ge} there exists a positive constant $C_0$ such that
\begin{equation}\label{pr9}
{\omega _{{D^*}}}\left( {0,\psi {{\left( {{F_\alpha }} \right)}^m}} \right) \ge {C_0}{e^{ - \pi {\lambda _{{D^*}}}\left( {{\gamma _m},\psi {{\left( {{F_\alpha }} \right)}^m}} \right)}}.
\end{equation}
Furthermore, by Lemma \ref{lem2} we infer that ${\gamma _m} \subset \overline D \left( {0,\frac{1}{2}} \right)$. So, taking the crosscuts ${\gamma _0} = \overline {{D^*}}  \cap \partial \mathbb{D}$ and ${\gamma _0}' = \overline {{D^*}}  \cap \partial D\left( {0,{e^{{{3\pi } \mathord{\left/
 {\vphantom {{3\pi } 2}} \right.
 \kern-\nulldelimiterspace} 2}}}} \right)$ of $D^*$ (see Fig. \ref{top1}, \ref{top2}) and applying Lemma \ref{lem3}  we obtain
 \begin{equation}\label{pr10}
 {\lambda _{{D^*}}}\left( {{\gamma _m},\psi {{\left( {{F_\alpha }} \right)}^m}} \right) \le 8.82{e^{ - 3\pi }} - 3 + {\lambda _{{D^*}}}\left( {{\gamma _m},{\gamma _0}'} \right) + {\lambda _{{D^*}}}\left( {{\gamma _0},\psi {{\left( {{F_\alpha }} \right)}^m}} \right).
 \end{equation}
 where ${\lambda _{{D^*}}}\left( {{\gamma _0},{\gamma _0}'} \right) = 3$ and ${\lambda _{{D^*}}}\left( {{\gamma _m},{\gamma _0}'} \right)$ is bounded from above by a positive constant $C_1$ for every $\alpha>0$ and every $m$ (see  \cite[p. 370-371]{Be} for a similar estimate). By relations (\ref{pr8}), (\ref{pr9}) and (\ref{pr10}), we get
 \begin{equation}\label{pr12}
 {\omega _\mathbb{D}}\left( {0,{F_\alpha }} \right) \ge {2^{n + 2}}{C_0}{e^{  \left( {3 - 8.82{e^{ - 3\pi }}} \right)\pi }}{e^{ -{C_1} \pi}}{e^{ - \pi {\lambda _{{D^*}}}\left( {{\gamma _0},\psi {{\left( {{F_\alpha }} \right)}^m}} \right)}}.
 \end{equation} 
 
\begin{figure}[H]
\begin{center}
\includegraphics[scale=0.5]{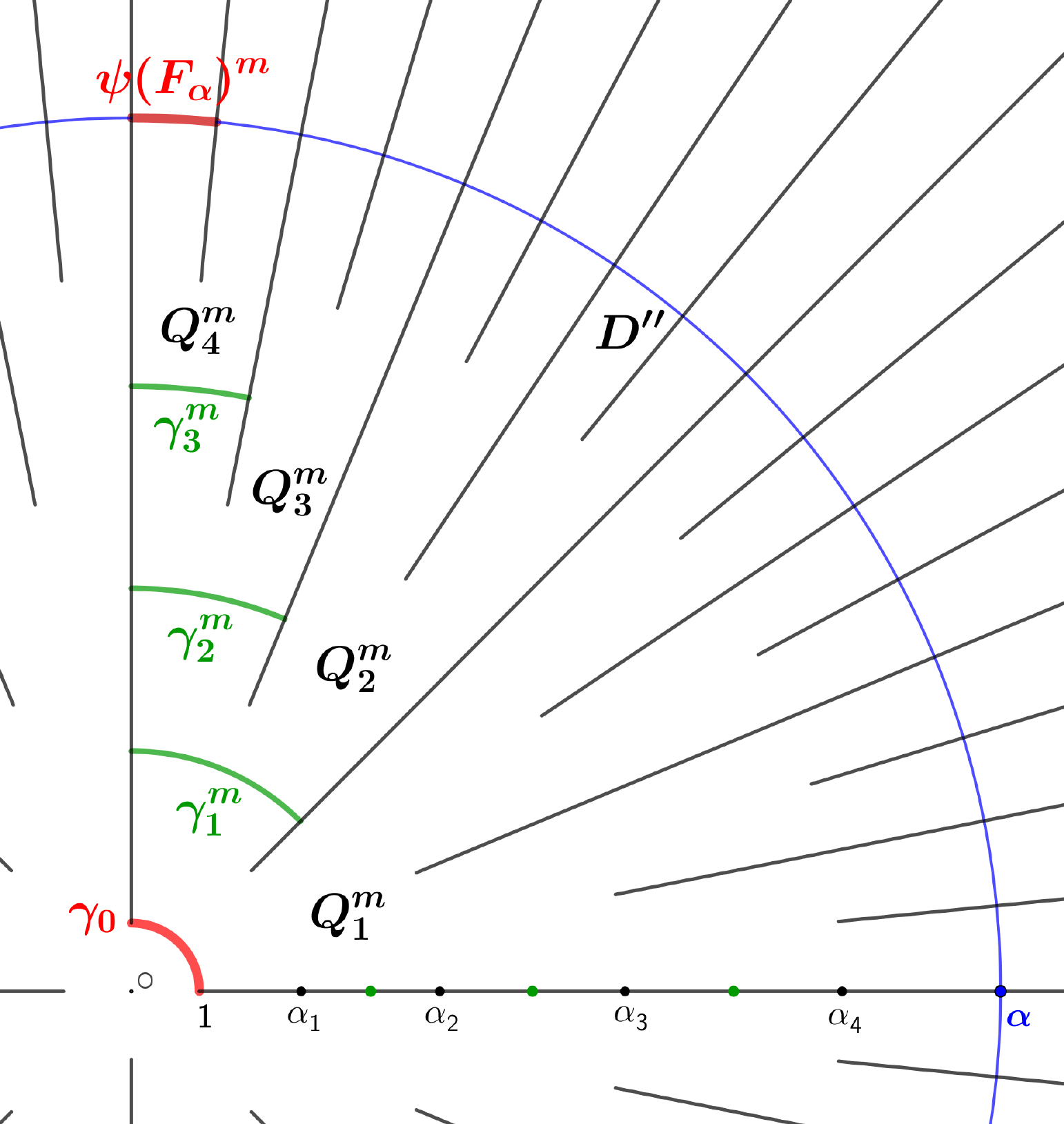}
\vspace*{0.3cm}
\caption{The simply connected domain $D''$ and the crosscuts $\gamma _1^m,\gamma _2^m,\gamma _3^m, \ldots ,\gamma _{n - 1}^m$ in case $\alpha  \in \left( {{\alpha _4},{\alpha _5}} \right)$.}
\label{deco}
\end{center}
\end{figure}

\textbf{Step 5:} Now we concentrate on ${e^{ - \pi {\lambda _{{{D^ * }}}}\left( {{\gamma _0},\psi {{\left( {{F_\alpha }} \right)}^m}} \right)}}$ or equivalently on ${e^{ - \pi {\lambda _{D''}}\left( {{\gamma _0},\psi {{\left( {{F_\alpha }} \right)}^m}} \right)}}$. First we take the crosscuts $\gamma _1^m,\gamma _2^m,\gamma _3^m, \ldots ,\gamma _{n - 1}^m$ of $D''$, where, for $j = 1,2,3, \ldots ,n - 1$, $\gamma _j^m$ is an arc of the circle with center at $0$ and radius equal to the midpoint of ${\left[ {{\alpha _j},{\alpha _{j + 1}}} \right]}$ as illustrated in Fig. \ref{deco}, and set
\[m\left( {{Q_1^m}} \right): = {\lambda _{D''}}\left( {{\gamma _0},{\gamma _1^m}} \right),\;m\left( {{Q_2^m}} \right): = {\lambda _{D''}}\left( {{\gamma _1^m},{\gamma _2^m}} \right), \ldots,\;m\left( {{Q_n^m}} \right): = {\lambda _{D''}}\left( {{\gamma _{n - 1}^m},{\psi \left( {F_\alpha } \right)^m}} \right).\]

\begin{figure}[H]
	\begin{center}
		\includegraphics[scale=0.33]{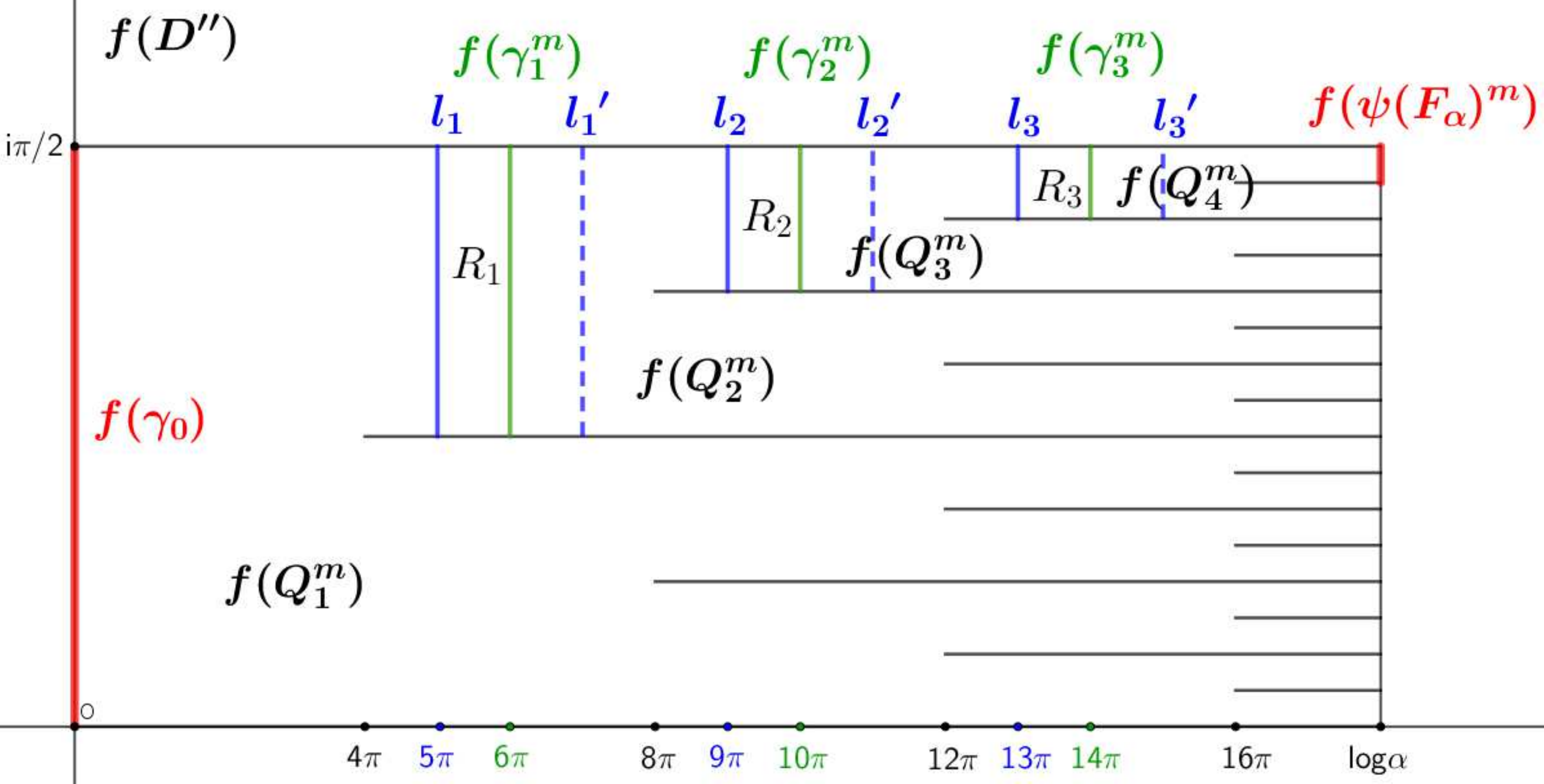}
		\vspace*{0.3cm}
		\caption{The image of $D''$ under the map $f\left( z \right) = \LOG z$ and its decomposition in case $\alpha  \in \left( {{\alpha _4},{\alpha _5}} \right)$.}
		\label{orth}
	\end{center}
\end{figure}

Applying the conformal map $f\left( z \right) = \LOG z$ on $D''$ we take $f\left( {D''} \right),f\left( {{\gamma _0}} \right),f\left( {\psi {{\left( {{F_\alpha }} \right)}^m}} \right)$ and $f\left( {\gamma _1^m} \right), \ldots ,f\left( {\gamma _{n - 1}^m} \right),f\left( {Q_1^m} \right), \ldots ,f\left( {Q_n^m} \right)$ illustrated in Fig. \ref{orth}. The conformal invariance of extremal length implies that for every $j = 1,2,3, \ldots ,n$,
\[{\lambda _{D''}}\left( {{\gamma _0},\psi {{\left( {{F_\alpha }} \right)}^m}} \right) = {\lambda _{f\left( {D''} \right)}}\left( {f\left( {{\gamma _0}} \right),f\left( {\psi {{\left( {{F_\alpha }} \right)}^m}} \right)} \right),\; m\left( {Q_j^m} \right) = m\left( {f\left( {Q_j^m} \right)} \right).\] 

This leads us to consider the crosscuts ${l_1},{l_1}',{l_2},{l_2}', \ldots ,{l_{n - 1}},{l_{n - 1}}'$ of $f\left( {D''} \right)$ so that, for $j = 1,2, \ldots ,n-1$, each of $l_j$ is a segment which lies on a vertical line passing through the midpoint of $\left[ {f\left( {{\alpha _j}} \right),\frac{{f\left( {{\alpha _j}} \right) + f\left( {{\alpha _{j + 1}}} \right)}}{2}} \right]=\left[ {4j \pi, 2\left( {2j+1} \right) \pi } \right]$ and ${l_j}'$ is the reflection of $l_j$ in $f\left( {\gamma _j^m} \right)$  (see Fig. \ref{orth}). Now let $R_j$ be the rectangle formed by ${\partial f\left( {D''} \right)}$  and ${l_j},f\left( {\gamma _j^m} \right)$ for $j = 1,2,3, \ldots ,n - 1$ as illustrated in Fig \ref{orth}. Since
\begin{eqnarray}
m\left( {{R_1}} \right) &=& \frac{{{\pi}}}{{{\pi  \mathord{\left/
				{\vphantom {\pi  4}} \right.
				\kern-\nulldelimiterspace} 4}}} = {2^2} > 3       \nonumber \\
m\left( {{R_2}} \right) &=& \frac{{{\pi}}}{{{\pi  \mathord{\left/
				{\vphantom {\pi  8}} \right.
				\kern-\nulldelimiterspace} 8}}} = {2^3} \nonumber \\
\vdots \nonumber\\
m\left( {{R_{n - 1}}} \right) &=& \frac{{{\pi }}}{{{\pi  \mathord{\left/
				{\vphantom {\pi  {{2^n}}}} \right.
				\kern-\nulldelimiterspace} {{2^n}}}}} = {2^{n}}, \nonumber
\end{eqnarray}
we can apply Lemma \ref{lem1} successively and obtain
\begin{eqnarray}
0 \le {\lambda _{D''}}\left( {{\gamma _0},\psi {{\left( {{F_\alpha }} \right)}^m}} \right) - \left( {m\left( {Q_1^m} \right) + m\left( {{{\left( {Q_1^m} \right)}^c}} \right)} \right) &\le& 26.46{e^{ - \pi m\left( {{R_1}} \right)}}      \nonumber \\
0 \le m\left( {{{\left( {Q_1^m} \right)}^c}} \right) - \left( {m\left( {Q_2^m} \right) + m\left( {{{\left( {Q_2^m} \right)}^c}} \right)} \right) &\le& 26.46{e^{ - \pi m\left( {{R_2}} \right)}} \nonumber \\
\vdots \nonumber \\
0 \le m\left( {{{\left( {Q_{n - 2}^m} \right)}^c}} \right) - \left( {m\left( {Q_{n - 1}^m} \right) + m\left( {{{\left( {Q_{n - 1}^m} \right)}^c}} \right)} \right) &\le& 26.46{e^{ - \pi m\left( {{R_{n - 1}}} \right)}},  \nonumber
\end{eqnarray}
where ${m\left( {{{\left( {Q_j^m} \right)}^c}} \right)}$ denotes the extremal length between $\gamma _j^m$ and $\psi {\left( {{F_\alpha }} \right)^m}$ in $D''$ for every $j = 1,2, \ldots ,n - 1$ and thus ${m\left( {{{\left( {Q_{n - 1}^m} \right)}^c}} \right) = m\left( {Q_n^m} \right)}$. Adding the inequalites above, we deduce that
\[ 0 \le {\lambda _{D''}}\left( {{\gamma _0},\psi {{\left( {{F_\alpha }} \right)}^m}} \right)\le 26.46\sum\limits_{j = 1}^{n - 1} {{e^{ - \pi m\left( {{R_j}} \right)}}}  + m\left( {Q_1^m} \right) + m\left( {Q_2^m} \right) +  \ldots  + m\left( {Q_n^m} \right),\]
where $m\left( {{R_j}} \right) = {2^{j+1}}$. So,
\[{e^{ - \pi {\lambda _{D''}}\left( {{\gamma _0},\psi {{\left( {{F_\alpha }} \right)}^m}} \right)}} \ge {e^{ - 26.46\pi \sum\limits_{j = 1}^{n - 1} {{e^{ - {2^{j + 1}}\pi }}} }}{e^{ - \pi \left( {m\left( {Q_1^m} \right) + m\left( {Q_2^m} \right) +  \ldots  + m\left( {Q_n^m} \right)} \right)}}.\]
Since the series ${\sum\limits_{j = 1}^{ + \infty } {{e^{ - {2^{j + 1}}\pi }}} }$ converges to a positive real number $l$, we obtain
\[{e^{ - \pi {\lambda _{D''}}\left( {{\gamma _0},\psi {{\left( {{F_\alpha }} \right)}^m}} \right)}} \ge {C_2}{e^{ - \pi \left( {m\left( {Q_1^m} \right) + m\left( {Q_2^m} \right) +  \ldots  + m\left( {Q_n^m} \right)} \right)}},\]
where $C_2 : ={e^{ - 26.46\pi l }}$. But due to the symmetry of $D$, we notice that the sum 
\[m\left( {Q_1^m} \right) + m\left( {Q_2^m} \right) +  \ldots  + m\left( {Q_n^m} \right)\]
is independent of $m$ and thus 
\[m\left( {Q_1^m} \right) + m\left( {Q_2^m} \right) +  \ldots  + m\left( {Q_n^m} \right) = m\left( {{Q_1}} \right) + m\left( {{Q_2}} \right) +  \ldots  + m\left( {{Q_n}} \right).\]
This leads us to the relation
\begin{equation}\label{pr13}
{e^{ - \pi {\lambda _{D''}}\left( {{\gamma _0},\psi {{\left( {{F_\alpha }} \right)}^m}} \right)}} \ge {C_2}{e^{ - \pi \left( {m\left( {{Q_1}} \right) + m\left( {{Q_2}} \right) +  \ldots  + m\left( {{Q_n}} \right)} \right)}}.
\end{equation}
Combining the relations (\ref{pr12}) and (\ref{pr13}) we finally obtain
\begin{equation}\label{pr14}
{\omega _\mathbb{D}}\left( {0,{F_\alpha }} \right) \ge {2^{n + 2}}{C_0}{C_2}{e^{  \left( {3 - 8.82{e^{ - 3\pi }}} \right)\pi }}{e^{ -{C_1} \pi}}{e^{ - \pi \left( {m\left( {{Q_1}} \right) + m\left( {{Q_2}} \right) +  \ldots  + m\left( {{Q_n}} \right)} \right)}}.
\end{equation}
\textbf{Step 6:} Now suppose there exists a positive constant $K$ such that for every $\alpha >0$,
\begin{equation}\label{pr15}
{\omega _\mathbb{D}}\left( {0,{F_\alpha }} \right) \le K{e^{ - {d_\mathbb{D}}\left( {0,{F_\alpha }} \right)}}.
\end{equation}
Setting
\[{C_3}: = \frac{{3{\pi ^2}}}{2}\frac{K}{{{{C_0}{C_2}}}}{e^{ {C_1} \pi}}{e^{-\left( {3 - 8.82{e^{ - 3\pi }}} \right)\pi }}\]
and using (\ref{pr7}), (\ref{pr14}) and (\ref{pr15}), we infer that
\[{2^{n + 2}}{e^{ - \pi \left( {m\left( {{Q_1}} \right) + m\left( {{Q_2}} \right) +  \ldots  + m\left( {{Q_n}} \right)} \right)}} \le {C_3}{e^{ - \pi \left( {m\left( {{Q_1}} \right) + m\left( {{Q_2}} \right) +  \ldots  + m\left( {{Q_n}} \right)} \right)}}\]
or equivalently
\begin{equation}\label{pr16}
{2^{n + 2}} \le {C_3}
\end{equation}
for every $n \in \mathbb{N}$. Finally, taking limits in (\ref{pr16}) as $n \to  + \infty $, that is $\alpha  \to  + \infty $, we obtain the contradiction
\[\mathop {\lim }\limits_{n \to  + \infty } {2^{n + 2}} \le {C_3} <  + \infty. \]
So, $\forall K>0$ $\exists \alpha$ such that
\[{\omega _\mathbb{D}}\left( {0,{F_\alpha }} \right) \ge K{e^{ - {d_\mathbb{D}}\left( {0,{F_\alpha }} \right)}}.\]
\qed

\section{The second example}\label{section5}

The main feature of $D$ which plays a central role in the proof of Theorem \ref{main} is that as $\alpha  \to  + \infty $, the number of components of $\psi \left( {{F_\alpha }} \right)$ tends to infinity. Next we prove that even if the number of components of $\psi \left( {{F_\alpha }} \right)$ is bounded from above by a positive constant for every $\alpha>0$, the answer to the Question \ref{quest} is still negative. To verify this, we need the following lemma whose proof is straightforward and thus is omitted. 

\begin{lemma}\label{geodesic}
	Let $\Gamma $ be the geodesic between two points $z_1,z_2 \in \partial{\mathbb{D}}$ in $\mathbb{D}$. Then
	\[{e^{ - {d_{\mathbb{D}}}\left( {0,\Gamma } \right)}} \le {\omega _{\mathbb{D}}}\left( {0,\Gamma } \right) \le \frac{4}{\pi }{e^{ - {d_{\mathbb{D}}}\left( {0,\Gamma } \right)}}.\]
\end{lemma}

\begin{theorem}\label{ex2} There exists an unbounded simply connected domain $D$ with the following properties: Let $\psi $ be a conformal map of $\mathbb{D}$ onto $D$. If ${F_\alpha } = \left\{ {z \in \mathbb{D}:\left| {\psi \left( z \right)} \right| = \alpha } \right\}$ for $\alpha >0$, then

\begin{enumerate}
\item $\psi \left( {{F_\alpha }} \right)$ is a connected set for every $\alpha>0$ and
\item $\forall K>0$ $\exists \alpha$ such that
\[{\omega _\mathbb{D}}\left( {0,{F_\alpha }} \right) \ge K{e^{ - {d_\mathbb{D}}\left( {0,{F_\alpha }} \right)}}.\]
\end{enumerate}

\end{theorem}

\begin{figure}[H] 
	\begin{center}
		\includegraphics[scale=0.3]{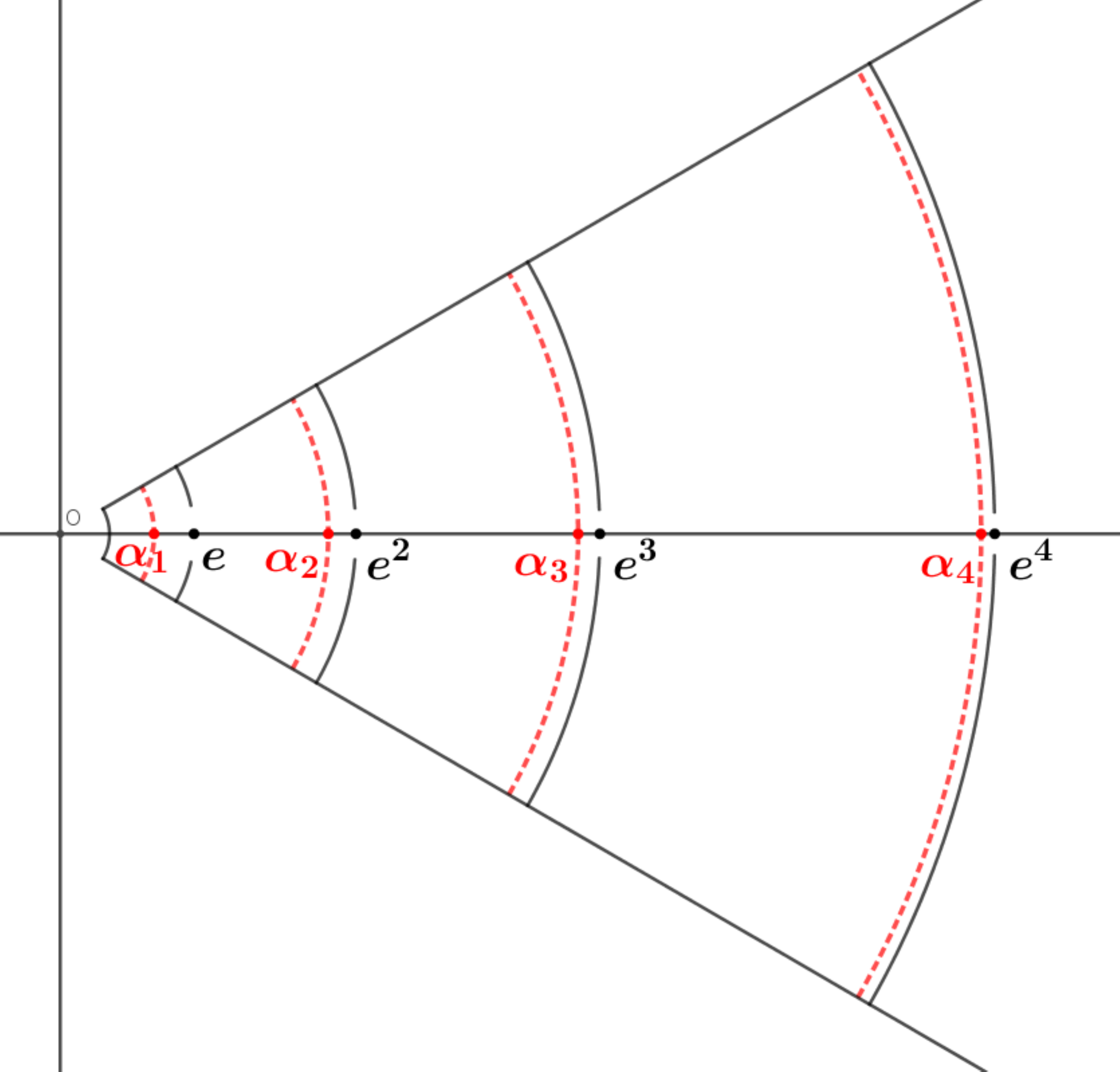}
		\caption{The simply connected domain $D$.}
		\label{ant}
	\end{center}
\end{figure}

\begin{figure}[H] 
	\begin{center}
		\includegraphics[scale=0.38]{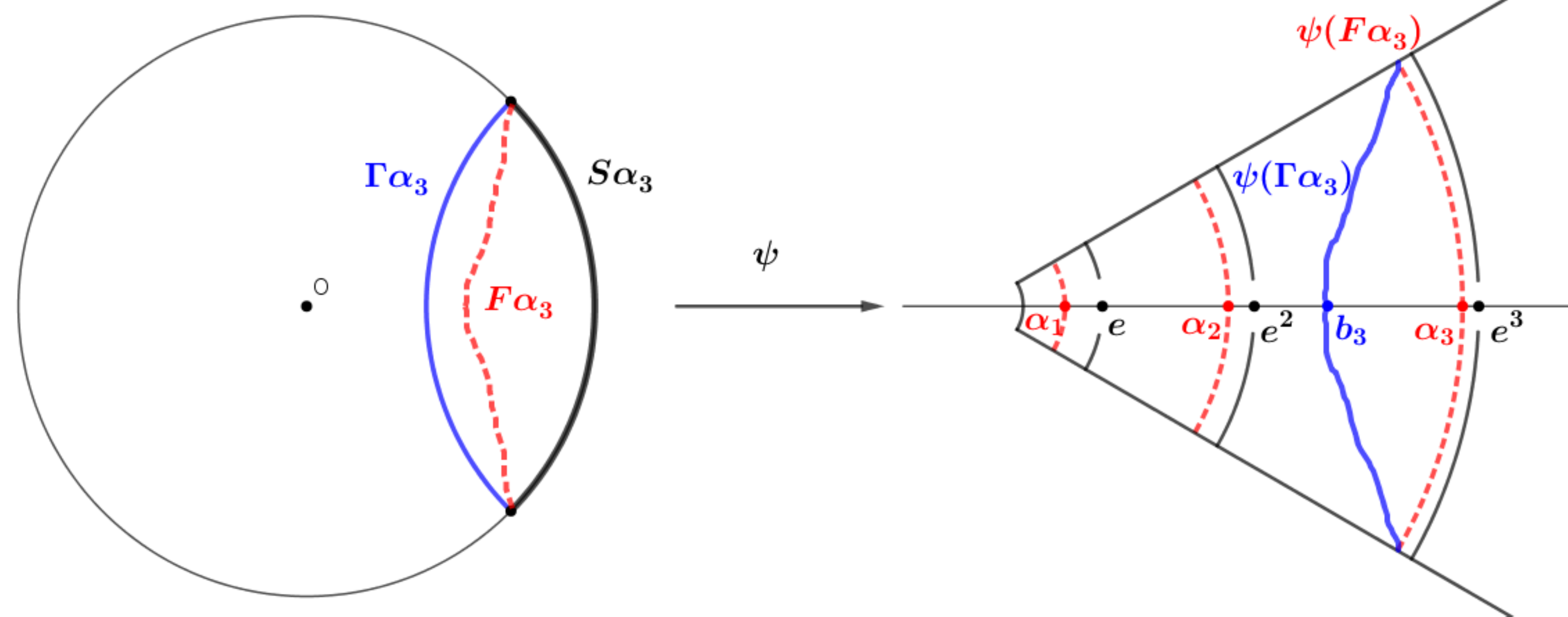}
		\caption{The curves $\Gamma _{{\alpha _n}},F_{\alpha _n},S_{\alpha _n}$ and their images under the map $\psi$ in case $n=3$.}
		\label{mapping}
	\end{center}
\end{figure}

\begin{figure}[H] 
	\begin{center}
		\includegraphics[scale=0.35]{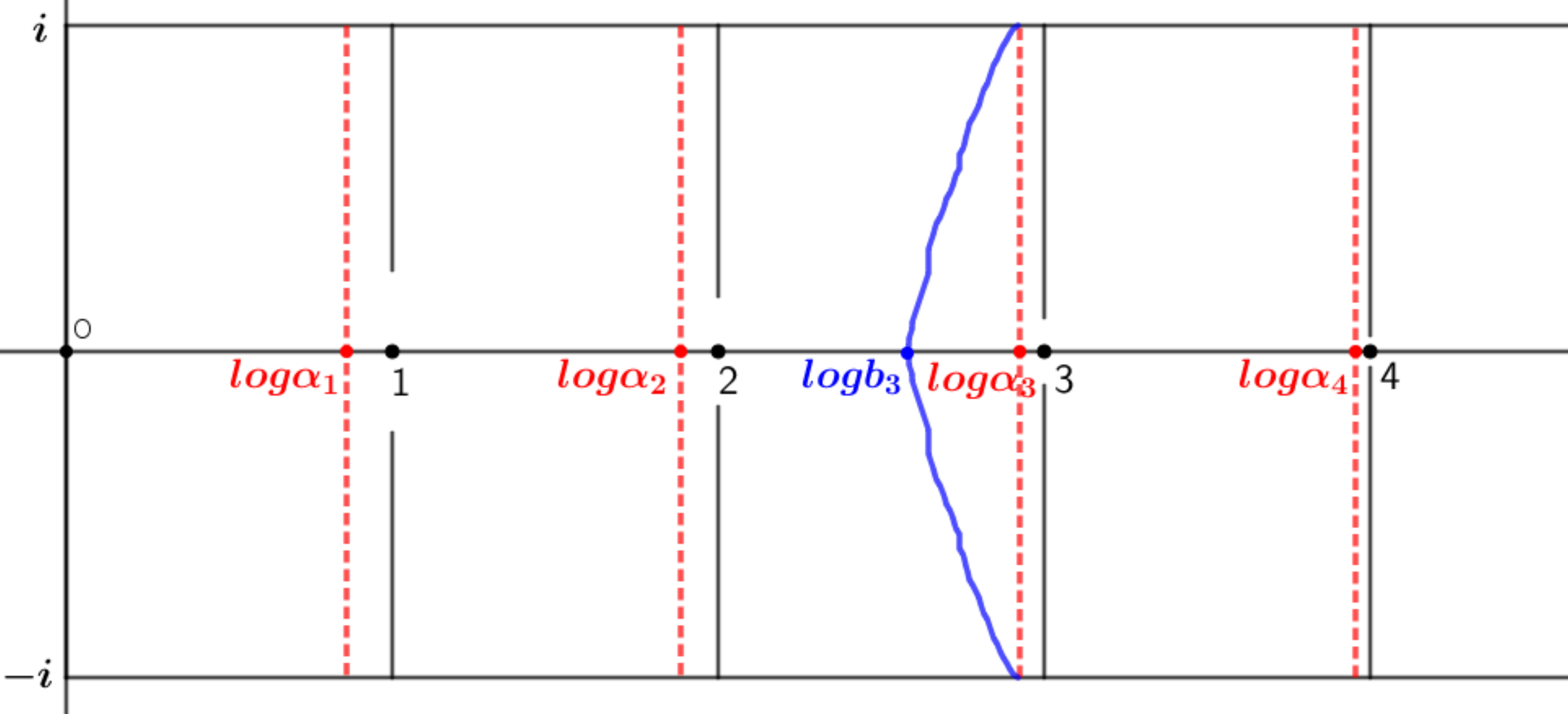}
		\caption{The simply connected domain $D'$.}
		\label{orty}
	\end{center}
\end{figure}

\begin{figure}[H] 
	\begin{center}
		\includegraphics[scale=0.3]{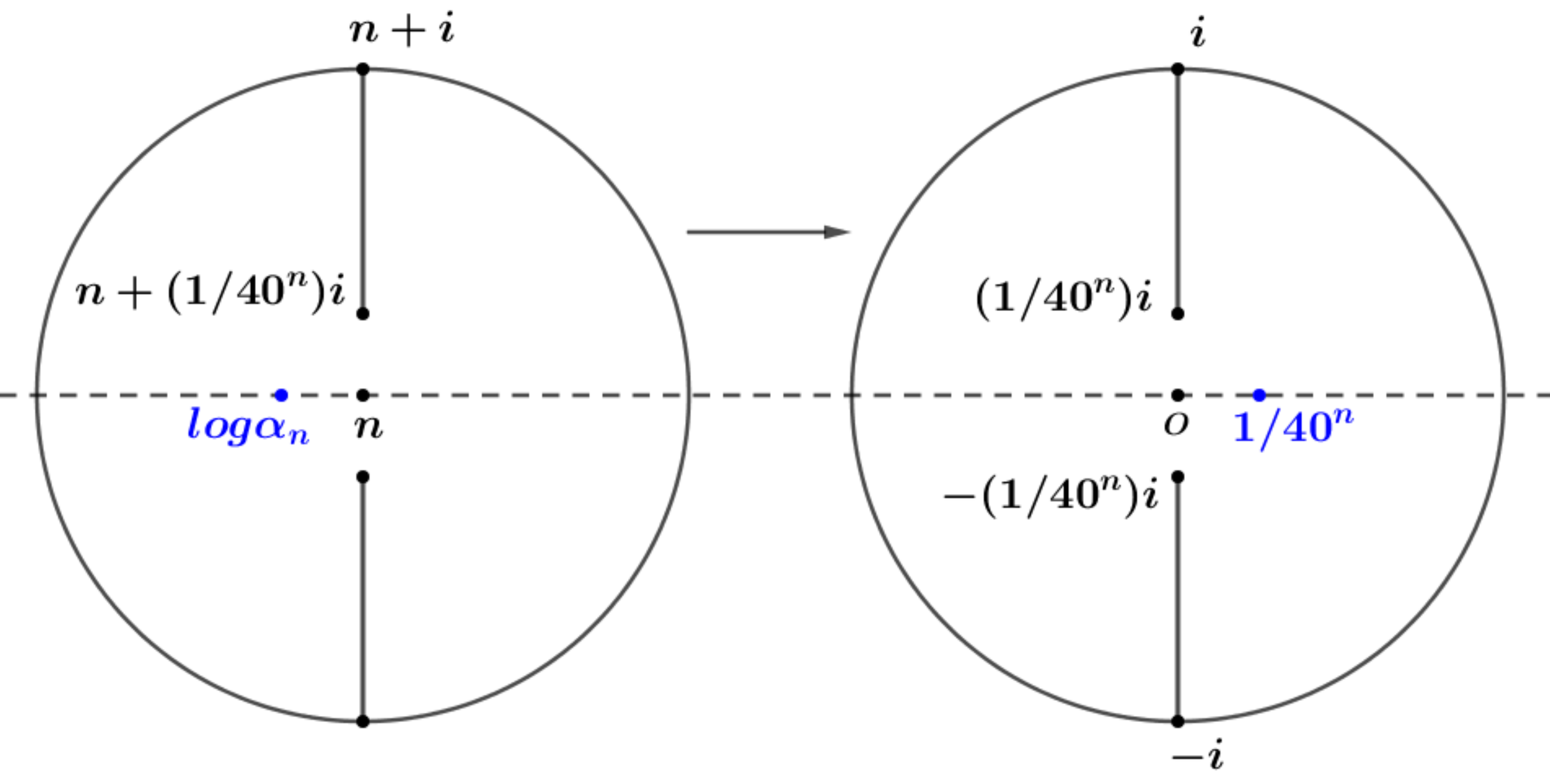}
		\caption{}
		\label{kik}
	\end{center}
\end{figure}

\proof Let $D$ be the simply connected domain of Fig. \ref{ant}, namely
\[D =   \left\{ {z \in \mathbb{C} \backslash \overline {\mathbb{D}} :\left| {\Arg{z}} \right| < 1} \right\}\backslash \bigcup\limits_{n = 1}^{ + \infty } {\left\{ {z \in \partial D\left( {0,{e^n}} \right):\frac{1}{{{40^n}}} \le \left| {\Arg{z}} \right| < 1} \right\}}.\]
The Riemann Mapping Theorem implies that there exists a conformal map $\psi $ from $\mathbb{D}$ onto $D$ such that $\psi \left( 0 \right) = {e^{\frac{1}{4}}}$. Let ${\alpha _n} = {e^{n - \frac{1}{{{{40}^n}}}}}$ for every $n \in \mathbb{N}$ and take the arcs ${\partial D\left( {0,{\alpha _n}} \right) \cap  \overline D }$ as illustrated in Fig. \ref{ant}. Now fix a number $n>1$. If ${\Gamma _{{\alpha _n}}}$ is the geodesic joining ${\psi ^{ - 1}}\left( {{\alpha _n}{e^i}} \right)$ to ${\psi ^{ - 1}}\left( {{\alpha _n}{e^{-i}}} \right)$ in $\mathbb{D}$ and $S _{{\alpha _n}}$ denotes the arc of $\partial \mathbb{D}$ between ${\psi ^{ - 1}}\left( {{\alpha _n}{e^{-i}}} \right)$ and ${\psi ^{ - 1}}\left( {{\alpha _n}{e^{i}}} \right)$ (see Fig. \ref{mapping}), then by Lemma \ref{geodesic} and \cite[p. 370]{Be} we get
\[{\omega _\mathbb{D}}\left( {0,{F_{{\alpha _n}}}} \right) \ge {\omega _\mathbb{D}}\left( {0,{S_{{\alpha _n}}}} \right) = \frac{1}{2}{\omega _\mathbb{D}}\left( {0,{\Gamma _{{\alpha _n}}}} \right) \ge \frac{1}{2}{e^{ - {d_\mathbb{D}}\left( {0,{\Gamma _{{\alpha _n}}}} \right)}} = \frac{1}{2}{e^{ - {d_D}\left( {{e^{\frac{1}{4}}},\psi \left( {{\Gamma _{{\alpha _n}}}} \right)} \right)}}.\]
Since $\psi$ preserves the geodesics and $D$ is symmetric with respect to the real axis, we deduce that ${d_D}\left( {{e^{\frac{1}{4}}},\psi \left( {{\Gamma _{{\alpha _n}}}} \right)} \right) = {d_D}\left( {{e^{\frac{1}{4}}},{b_n}} \right)$, where ${{b_n} \in \left( {{e^{n - 1}},{e^n}} \right)}$. So,
\begin{equation}\label{sxesi1}
{\omega _\mathbb{D}}\left( {0,{F_{{\alpha _n}}}} \right) \ge \frac{1}{2}{e^{ - {d_D}\left( {{e^{\frac{1}{4}}},{b_n}} \right)}}.
\end{equation}
Notice that if ${g_D}\left( {{e^{\frac{1}{4}}},z} \right)$ denotes the Green's function for $D$ (see \cite[p. 41-43]{Gar}, \cite[p. 106-115]{Ra}), then
\[{d_D}\left( {{e^{\frac{1}{4}}},z} \right) = \log \frac{{1 + {e^{ - {g_D}\left( {{e^{\frac{1}{4}}},z} \right)}}}}{{1 - {e^{ - {g_D}\left( {{e^{\frac{1}{4}}},z} \right)}}}}\]
(see \cite[p. 12-13]{Bea} and \cite[p. 106]{Ra}). Consider the conformal map $h\left( z \right) = \LOG z$ of $D$ onto $ D':= h\left( D \right)$ (see Fig. \ref{orty}). For every $\alpha {'_n} \in \psi \left( {{F_{{\alpha _n}}}} \right)\backslash \left\{ {{\alpha _n}} \right\}$, we infer, by a symmetrization result, that
\[{{g_D}\left( {{e^{\frac{1}{4}}},{\alpha _n}} \right)={g_{D'}}\left( {\frac{1}{4}},\log{{\alpha _n}} \right) \ge {g_{D'}}\left( {{\frac{1}{4}},\log {\alpha {'_n}}} \right)={g_D}\left( {{e^{\frac{1}{4}}},\alpha {'_n}} \right)},\]
(see Lemma 9.4 \cite[p. 659]{Hay}). Since
\[f\left( x \right) = \log \frac{{1 + {e^{ - x}}}}{{1 - {e^{ - x}}}}\]
is a decreasing function on $\left( {0, + \infty } \right)$, we have that ${d_D}\left( {{e^{\frac{1}{4}}},\psi \left( {{F_{{\alpha _n}}}} \right)} \right) = {d_D}\left( {{e^{\frac{1}{4}}},{\alpha _n}} \right)$. Thus,
\begin{equation}\label{sxesi2}
{e^{ - {d_\mathbb{D}}\left( {0,{F_{{\alpha _n}}}} \right)}} = {e^{ - {d_D}\left( {{e^{\frac{1}{4}}},\psi \left( {{F_{{\alpha _n}}}} \right)} \right)}} = {e^{ - {d_D}\left( {{e^{\frac{1}{4}}},{\alpha _n}} \right)}}.
\end{equation}
Next we prove that $b_n< \alpha_n$. Since $b_n$ lies on the geodesic $\psi \left( {{{\Gamma}_{{\alpha _n}}}} \right)$, we have
that ${\omega _D}\left( {{b_n},\psi \left( {{S_{{\alpha _n}}}} \right)} \right) = \frac{1}{2}$. Denoting by $D\left( {n,1} \right)=\left\{ {z \in \mathbb{C}:\left| {z - n} \right| < 1} \right\}$, the conformal invariance of harmonic measure implies that
\[{\omega _D}\left( {{\alpha _n},\psi \left( {{S_{{\alpha _n}}}} \right)} \right) = {\omega _{D'}}\left( {\log {\alpha _n},h \circ \psi \left( {{S_{{\alpha _n}}}} \right)} \right) \ge {\omega _{D\left( {n,1} \right)}}\left( {\log {\alpha _n},\left\{ {n + iy:\frac{1}{{{{40}^n}}} \le \left| y \right| \le 1} \right\}} \right),\]
where the last term, through a translation and a rotation (see Fig. \ref{kik}), can be expressed by
\[{\omega _{D\left( {n,1} \right)}}\left( {\log {\alpha _n},\left\{ {n + iy:\frac{1}{{{{40}^n}}} \le \left| y \right| \le 1} \right\}} \right) = {\omega _\mathbb{D}}\left( {\frac{1}{{{{40}^n}}},\left\{ {iy:\frac{1}{{{{40}^n}}} \le \left| y \right| \le 1} \right\}} \right).\]
Applying Beurling-Nevanlinna projection theorem \cite[p. 43]{Ahl}, we get
\[{\omega _\mathbb{D}}\left( {\frac{1}{{{{40}^n}}},\left\{ {iy:\frac{1}{{{{40}^n}}} \le \left| y \right| \le 1} \right\}} \right) \ge {\omega _\mathbb{D}}\left( {\frac{1}{{{{40}^n}}},\left[ { - 1, - \frac{1}{{{{40}^n}}}} \right]} \right) \ge {\omega _\mathbb{D}}\left( {\frac{1}{{40}},\left[ { - 1, - \frac{1}{{40}}} \right]} \right),\]
where
\[{\omega _\mathbb{D}}\left( {\frac{1}{{40}},\left[ { - 1, - \frac{1}{{40}}} \right]} \right) = \frac{2}{\pi }\arcsin {\left( {\frac{{1 - \frac{1}{{40}}}}{{1 + \frac{1}{{40}}}}} \right)^2}=0.719987303 > 0.7\]
Therefore,
\[{\omega _D}\left( {{\alpha _n},\psi \left( {{S_{{\alpha _n}}}} \right)} \right) > 0.7 > 0.5 = {\omega _D}\left( {{b_n},\psi \left( {{S_{{\alpha _n}}}} \right)} \right)\]
which implies that $b_n< \alpha_n$ and thus
\begin{equation}\label{isotita}
{e^{ - {d_D}\left( {{e^{\frac{1}{4}}},{\alpha _n}} \right)}} = {e^{ - {d_D}\left( {{e^{\frac{1}{4}}},{b_n}} \right)}}{e^{ - {d_D}\left( {{b_n},{\alpha _n}} \right)}}.
\end{equation}
Now suppose there exists a positive constant $K$ such that for every $\alpha >0$,
\begin{equation}\label{sxesi3}
{\omega _\mathbb{D}}\left( {0,{F_\alpha }} \right) \le K{e^{ - {d_\mathbb{D}}\left( {0,{F_\alpha }} \right)}}.
\end{equation}
Combining the relations (\ref{sxesi1}), (\ref{sxesi2}), (\ref{isotita}) and (\ref{sxesi3}), we infer that
\[\frac{1}{2}{e^{ - {d_D}\left( {{e^{\frac{1}{4}}},{b_n}} \right)}} \le {\omega _\mathbb{D}}\left( {0,{F_{{\alpha _n}}}} \right) \le K{e^{ - {d_\mathbb{D}}\left( {0,{F_{{\alpha _n}}}} \right)}} = K{e^{ - {d_D}\left( {{e^{\frac{1}{4}}},{b_n}} \right)}}{e^{ - {d_D}\left( {{b_n},{\alpha _n}} \right)}}\]
or
\begin{equation}\label{atopo}
{e^{  {d_D}\left( {{b_n},{\alpha _n}} \right)}} \le 2K
\end{equation}
for every $n>1$. However, using the quasi-hyperbolic distance defined in Section \ref{section2}, we get
\begin{eqnarray}\label{anisot}
{d_D}\left( {{b_n},{\alpha _n}} \right)&=& {d_{D'}}\left( {\log {b_n},\log{\alpha _n}} \right) \ge \frac{1}{2}{\delta _{D'}}\left( {\log {b_n},\log{\alpha _n}} \right) = \frac{1}{2}\int_{{{\mathop{\log {b _n}}} }}^{\log {\alpha _n}} {\frac{{dx}}{{d\left( {x,\partial D'} \right)}}}        \nonumber \\
&\ge&\frac{1}{2}\int_{\log {b_n}}^{\log{\alpha _n}} {\frac{{dx}}{{\sqrt {{{\left( {\frac{1}{{{40^{n}}}}} \right)}^2} + {{\left( {n - x} \right)}^2}} }}}  = \left. { - \frac{1}{2}\arcsinh \left( {{40^{n}}\left( {n - x} \right)} \right)} \right|_{\log {b_n}}^{\log{\alpha _n}} \nonumber \\
&=& - \frac{1}{2}\arcsinh \left( {1} \right) + \frac{1}{2}\arcsinh \left( {{40^{n}}\left( {n-\log b_n} \right)} \right) \nonumber \\
&\ge&  - \frac{1}{2}\arcsinh \left( {{1}} \right) + \frac{1}{2}\arcsinh \left( {{40^{n}}k} \right), 
\end{eqnarray}
where $k$ is a positive constant whose existence comes from the fact that for every $n>1$,
\[n - \log {b_n} \ge \log {\alpha _n} - \log {b_n}\] 
and 
\[{\omega _{D'}}\left( {\log {b_n},h \circ \psi \left( {{S_{{\alpha _n}}}} \right)} \right) = \frac{1}{2} < 0.7 < {\omega _{D'}}\left( {\log {\alpha _n},h \circ \psi \left( {{S_{{\alpha _n}}}} \right)} \right),\]
as we proved before. Finally, taking limits in (\ref{anisot}) as $n \to  + \infty $, we obtain the contradiction to the relation (\ref{atopo}),
\[\mathop {\lim }\limits_{n \to  + \infty } {d_D}\left( {{b_n},{\alpha _n}} \right) =  + \infty. \]
So, $\forall K>0$ $\exists \alpha$ such that
\[{\omega _\mathbb{D}}\left( {0,{F_\alpha }} \right) \ge K{e^{ - {d_\mathbb{D}}\left( {0,{F_\alpha }} \right)}}.\]
\qed

\section{Proof of Theorem \ref{kyrio}}\label{section6}

\begin{figure}[H] 
	\begin{center}
		\includegraphics[scale=0.28]{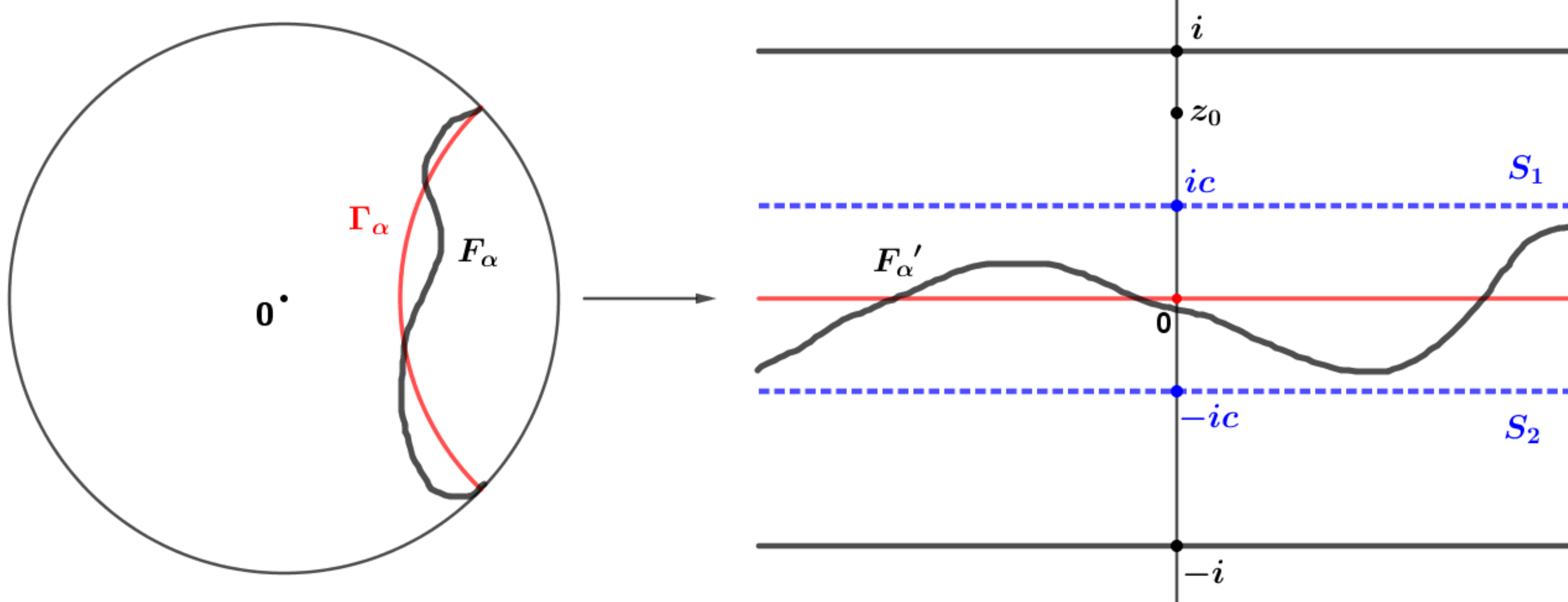}
		\caption{The conformal mapping of $\mathbb{D}$ onto the strip $S$.}
		\label{theor}
	\end{center}
\end{figure}

\proof Because of the assumption (1) and the additivity of harmonic measure we may assume that $N\left( \alpha  \right) =1$. We map conformally $\mathbb{D}$ onto the strip $S = \left\{ {z \in \mathbb{C}:\left| {{\mathop{\IM z}\nolimits} } \right| < 1} \right\}$ so that 
\[0 \mapsto {z_0} \in i \mathbb{R}^+,\;{\Gamma _\alpha } \mapsto \mathbb{R}.\]
Let ${F_\alpha }'$ be the image of ${F_\alpha }$. By assumption (2), there exists a positive constant $c=c \left( c_2  \right)$ such that $c < \left| {{z_0}} \right|$ and ${F_\alpha }' \subset \left\{ {z \in S:\left|{{\mathop{\IM z}\nolimits} } \right| < c} \right\}$ for every $\alpha>0$ (see Fig. \ref{theor}).
Set ${S_1} = \left\{ {z \in S:{\mathop{\IM z}\nolimits}  = c} \right\}$ and ${S_2} = \left\{ {z \in S:{\mathop{\IM z}\nolimits}  = -c} \right\}$. Then we have
\begin{equation}\label{rel1}
{\omega _{\mathbb{D}}}\left( {0,{F_\alpha }} \right) = {\omega _S}\left( {{z_0},{F_\alpha }'} \right) \le {\omega _S}\left( {{z_0},{S_1}} \right).
\end{equation}
Notice that by symmetry, for every $z \in S_1$,
\[{\omega _S}\left( {z,{S_2}} \right) = {\omega _S}\left( {ic,{S_2}} \right) = \frac{{1 - c}}{{1 + c}},\]
where the second equality comes from \cite[p. 100]{Ra}. Therefore, the strong Markov property for harmonic measure (see \cite[p. 282]{Bet}) implies that
\[{\omega _S}\left( {{z_0},{S_2}} \right) = \int_{{S_1}} {{\omega _S}\left( {z,{S_2}} \right){\omega _S}\left( {{z_0},dz} \right)}  = {\omega _S}\left( {{z_0},{S_1}} \right){\omega _S}\left( {ic,{S_2}} \right)\]
or
\begin{equation}\label{rel2}
{\omega _S}\left( {{z_0},{S_1}} \right) = \frac{{1 + c}}{{1 - c}}\;{\omega _S}\left( {{z_0},{S_2}} \right).
\end{equation}
Combining the relations (\ref{rel1}) and (\ref{rel2}), we get
\begin{equation}\label{rel3}
{\omega _{\mathbb{D}}}\left( {0,{F_\alpha }} \right) \le \frac{{1 + c}}{{1 - c}}\;{\omega _S}\left( {{z_0},{S_2}} \right) \le \frac{{1 + c}}{{1 - c}}\;{\omega _S}\left( {{z_0},\mathbb{R}} \right).
\end{equation}
Conformal invariance and Lemma \ref{geodesic} imply that 
\[{\omega _S}\left( {{z_0},\mathbb{R}} \right) \le \frac{4}{\pi }{e^{ - {d_S}\left( {{z_0},\mathbb{R}} \right)}}\]
which in conjunction with (\ref{rel3}) leads to
\begin{equation}\label{rel4}
{\omega _{\mathbb{D}}}\left( {0,{F_\alpha }} \right) \le \frac{4}{\pi } \frac{{1 + c}}{{1 - c}}\;{e^{ - {d_S}\left( {{z_0},\mathbb{R}} \right)}} \le \frac{4}{\pi }  \frac{{1 + c}}{{1 - c}}\;{e^{ - {d_S}\left( {{z_0},{S_1}} \right)}}.
\end{equation}
But
\begin{eqnarray}\label{rel5}
{d_S}\left( {{z_0},{S_2}} \right) &=& {d_S}\left( {{z_0}, - ic} \right) = {d_S}\left( {{z_0},ic} \right) + {d_S}\left( {ic, - ic} \right)       \nonumber \\
&=& {d_S}\left( {{z_0},S_1} \right) + {d_S}\left( {-ic,ic} \right),
\end{eqnarray}
where by \cite[p. 31]{Bea},
\begin{equation}\label{rel6}
{d_S}\left( { - ic,ic} \right) = \int_{ - \frac{{c\pi }}{2}}^{\frac{{c\pi }}{2}} {\frac{{dt}}{{\cos t}}}  = \log \frac{{1 + \sin \left( {{{c\pi } \mathord{\left/
 {\vphantom {{c\pi } 2}} \right.
 \kern-\nulldelimiterspace} 2}} \right)}}{{1 - \sin \left( {{{c\pi } \mathord{\left/
 {\vphantom {{c\pi } 2}} \right.
 \kern-\nulldelimiterspace} 2}} \right)}}.
\end{equation}
Combining the relations (\ref{rel5}) and (\ref{rel6}), we infer that
\[ - {d_S}\left( {{z_0},{S_1}} \right) =  - {d_S}\left( {{z_0},{S_2}} \right) + \log \frac{{1 + \sin \left( {{{c\pi } \mathord{\left/
 {\vphantom {{c\pi } 2}} \right.
 \kern-\nulldelimiterspace} 2}} \right)}}{{1 - \sin \left( {{{c\pi } \mathord{\left/
 {\vphantom {{c\pi } 2}} \right.
 \kern-\nulldelimiterspace} 2}} \right)}} \le  - {d_S}\left( {{z_0},{F_\alpha }'} \right) + \log \frac{{1 + \sin \left( {{{c\pi } \mathord{\left/
 {\vphantom {{c\pi } 2}} \right.
 \kern-\nulldelimiterspace} 2}} \right)}}{{1 - \sin \left( {{{c\pi } \mathord{\left/
 {\vphantom {{c\pi } 2}} \right.
 \kern-\nulldelimiterspace} 2}} \right)}}.\]
This together with (\ref{rel4}) give
\[{\omega _{\mathbb{D}}}\left( {0,{F_\alpha }} \right) \le \frac{4}{\pi }  \frac{{1 + c}}{{1 - c}}  \frac{{1 + \sin \left( {{{c\pi } \mathord{\left/
 {\vphantom {{c\pi } 2}} \right.
 \kern-\nulldelimiterspace} 2}} \right)}}{{1 - \sin \left( {{{c\pi } \mathord{\left/
 {\vphantom {{c\pi } 2}} \right.
 \kern-\nulldelimiterspace} 2}} \right)}}\;{e^{ - {d_S}\left( {{z_0},{F_\alpha }'} \right)}} = \frac{4}{\pi }  \frac{{1 + c}}{{1 - c}} \frac{{1 + \sin \left( {{{c\pi } \mathord{\left/
 {\vphantom {{c\pi } 2}} \right.
 \kern-\nulldelimiterspace} 2}} \right)}}{{1 - \sin \left( {{{c\pi } \mathord{\left/
 {\vphantom {{c\pi } 2}} \right.
 \kern-\nulldelimiterspace} 2}} \right)}}\;{e^{ - {d_{\mathbb{D}}}\left( {0,{F_\alpha }} \right)}}.\]
Thus, setting $K: = \frac{4}{\pi } \frac{{1 + c}}{{1 - c}} \frac{{1 + \sin \left( {{{c\pi } \mathord{\left/
  {\vphantom {{c\pi } 2}} \right.
  \kern-\nulldelimiterspace} 2}} \right)}}{{1 - \sin \left( {{{c\pi } \mathord{\left/
  {\vphantom {{c\pi } 2}} \right.
  \kern-\nulldelimiterspace} 2}} \right)}}$, we finally get that for every $\alpha>0$,
\[{\omega _\mathbb{D}}\left( {0,{F_\alpha }} \right) \le K{e^{ - {d_\mathbb{D}}\left( {0,{F_\alpha }} \right)}}.\]
\qed

\end{document}